\documentclass[review,11pt]{llncs}
\usepackage{graphicx}

\usepackage{lineno}

\usepackage{color}

\modulolinenumbers[5]

\usepackage{amssymb}
\usepackage{amsmath}
\usepackage{tikz}
\usepackage{subfig}
\usepackage{caption}
\usepackage{float}
\usepackage{mathrsfs}
\usepackage{eufrak}
\usepackage{enumerate}

\usepackage[colorlinks=true,breaklinks=true,bookmarks=true,urlcolor=blue,
     citecolor=blue,linkcolor=blue,bookmarksopen=false,draft=false]{hyperref}

\makeatletter
\tikzset{my loop/.style =  {to path={
  \pgfextra{}
  [looseness=8,min distance=7mm]
  \tikz@to@curve@path},font=\sffamily\small
  }}  
\makeatletter

\usepackage{setspace}

\newtheorem{theorem*}{Theorem}
\newcommand{\R}{\mathbb{R}}
\newcommand{\supp}{\text{supp}}
\newcommand{\spc}{\text{ }}

\usepackage{algorithm}

\usepackage{natbib}

\definecolor{OliveGreen}{cmyk}{0.6,0,0.95,0.40} 

\begin{document}



\title{An LP-based, Strongly-Polynomial 2-Approximation\\ Algorithm for Sparse Wasserstein Barycenters}


\author{Steffen Borgwardt}
\institute{\email{\href{mailto:steffen.borgwardt@ucdenver.edu}{steffen.borgwardt@ucdenver.edu}} 
University of Colorado Denver
}






\maketitle

\begin{abstract}
Discrete Wasserstein barycenters correspond to optimal solutions of transportation problems for a set of probability measures with finite support. Discrete barycenters are measures with finite support themselves and exhibit two favorable properties: there always exists one with a provably sparse support, and any optimal transport to the input measures is non-mass splitting. 

It is open whether a discrete barycenter can be computed in polynomial time. It is possible to find an exact barycenter through linear programming, but these programs may scale exponentially. In this paper, we prove that there is a strongly-polynomial $2$-approximation algorithm based on linear programming. First, we show that an exact computation over the union of supports of the input measures gives a tight $2$-approximation. This computation can be done through a linear program with setup and solution in strongly-polynomial time. The resulting measure is sparse, but an optimal transport may split mass. We then devise a second, strongly-polynomial algorithm to improve this measure to one with a non-mass splitting transport of lower cost. The key step is an update of the possible support set to resolve mass split.

Finally, we devise an iterative scheme that alternates between these two algorithms. The algorithm terminates with a $2$-approximation that has both a sparse support and an associated non-mass splitting optimal transport. We conclude with some sample computations and an analysis of the scaling of our algorithms, exhibiting vast improvements in running time over exact LP-based computations and low practical errors.

\end{abstract}




\section{Introduction}
Transportation problems for several marginals arise in applications ranging from finance and economics \citep{bhp-13,ght-14,rpdb-12,p-14} to physics \citep{bpg-12,cfk-13}, economics \citep{ce-10,cmn-10}, statistics \citep{bk-17,bll-11,mmh-11}, and data analytics \citep{bcmm-16,ggr-18}. The so-called {\em Wasserstein barycenters} correspond to optimal solutions to these problems, and have seen much recent attention. Barycenters are intimately connected to Fr\'echet means in Euclidean space \citep{mtbmmh-15,ty-05,tmmh-14,zp-17}, which is one of the origins of this field of research and the reason why statistical and probability notation is commonly used.

Given probability measures $P_1,\ldots, P_N$ on $\Bbb R^d$ and a weight vector $\lambda=(\lambda_1,\dots,\lambda_N) \in \R^N_{> 0}$ with $\sum_{i=1}^N \lambda_i=1$, a ($\lambda$-weighted) Wasserstein barycenter is a probability measure $\bar P$ on $\Bbb R^d$ which satisfies
\begin{equation}
\label{barycenter}
\phi(\bar P):= \sum_{i=1}^N \lambda_i W_2(\bar{P},P_i)^2=\inf_{P\in\mathcal P^2(\Bbb R^d)}\sum_{i=1}^N \lambda_i W_2( P, P_i)^2,
\end{equation}
where $W_2$ is the quadratic Wasserstein distance and  $\mathcal P^2(\Bbb R^d)$ is the set of all probability measures on $\Bbb R^d$ with finite second moments. We recommend the monographs \cite{v-03,v-09}, and the more recent \cite{pz-19,pc-19}, for a review of the Wasserstein distance and an overview of the literature on optimal transport problems. 

\subsection{Exact, approximate and heuristic algorithms}\label{sec:exact}

Exact barycenter computations for a set of continuous measures are intractable outside of some special cases, in particular because an evaluation of the Wasserstein distance between two given continuous measures itself already is challenging. Because of this, the literature uses several types of simplifications to facilitate practical computations. The arguably most important one is an input of discrete data:
 
In many applications, data is given as a set of {\em discrete (probability) measures} $P_1,\ldots, P_N$ having finite support in $\mathbb{R}^d$. A {\em discrete Wasserstein barycenter} is a probability measure $\bar{P}$ which satisfies Eq.~(\ref{barycenter}) for such measures. Discrete probability measures arise naturally in applications in operations research. The finite support often corresponds to a set of geographical locations (customers, facilities, service providers) and the measures represent data that varies over time. In computer graphics or image science, discrete probability measures are supported on a grid or arise through a discretization of the underlying space.

In \cite{abm-16}, some theoretical results were developed for these discrete barycenters. They mirror the continuous case, established in \cite{ac-11}, with a few notable exceptions. First, unlike in the continuous case, there may exist several discrete barycenters for the same set of measures. All of them have finite support and there always exists a discrete barycenter with provably sparse support. Analogously to the continuous case, a discrete barycenter always has a non-mass splitting optimal transport to each discrete marginal, i.e., each barycenter support point transports its whole mass to just a single support point for each measure. These results were proven for uniform $\lambda_i=\frac{1}{n}$ in \cite{abm-16}, but are readily transferred to general $\lambda_i$ (\cite{m-16}). 

Both sparsity and non-mass split are crucial to applications. Sparsity is desirable in many applications of operations research such as facility location. Non-mass split is often imposed by physical limitations of applications, such as the design of deformable templates \citep{bll-11,jzd-98,ty-05}. For example, in metal shaping, sheets of metal have to be pressed into a collection of different shapes. Each of these shapes is modeled as a measure. A `mean deformation' (barycenter) is a best shape for the initial sheet of metal with respect to the energy required to mold (transport) it into all required shapes. Only because of the existence of a non-mass splitting transport, the mean deformation can indeed be transformed into each shape through only bending and stretching. See \cite{bll-11} for more details on such applications. We are interested in the computation of a discrete barycenter that exhibits these beneficial properties. 

It is open whether the computation of an exact discrete barycenter can be done in polynomial time. It is well-known that linear programming can be used to approximate or solve optimal transport problems \citep{sld-18,ywwl-17}. Most importantly, exact discrete barycenters can be computed through linear programming \citep{abm-16,bs-18,coo-15}. However, these programs may scale exponentially in the number of measures $N$ (see Section \ref{stateoftheart2}) and thus have not been widely considered for practical use. Much larger and faster practical computations are possible through various heuristics that fundamentally differ from LP-based approaches. Most of these algorithms are based on simplifications to the Wasserstein distance to obtain an easier objective function. 

The arguably most popular tools are entropic regularization techniques, which are used to make the objective function smooth and strictly convex \citep{c-13}. In recent years there has been significant progress on these techniques and they led to a flurry of competitive algorithms for good approximations of barycenters in practice. For example, regularization is used in a well-behaved implementation of a gradient descent algorithm that uses information from both smoothed primal and dual optimal transport formulations \citep{cd-14}. Further, regularization not only greatly simplifies the underlying optimal transport problem itself, but the regularized barycenter problem then also allows the efficient computation of iterative Bregman projections \citep{bccnp-14}. Iterative Bregman methods have proven to be a competitive approach for large-scale computations \citep{ylst-19,ywwl-17}; their number of variables scales roughly linearly in the number of marginals. The entropy regularized Wasserstein distance converges towards the actual Wasserstein distance in $O(\frac{1}{w})$, where $w$ is the entropic regularization factor, and the non-regularized transport cost computed with a regularized transport plan converges towards the Wasserstein distance in $O(e^{-w})$ \citep{bccnp-14,cdps-17,lrpc-18}. The factor $w$ usually is chosen empirically.

The great scalability of regularization-based methods, see for example \cite{sgpcbndg-15}, comes at the cost of a few drawbacks. First, they typically require a fixed support over which a barycenter approximation is to be computed. For grid-structured data, approaches in the literature usually just specify the underlying grid as the support set. One of the main results in this paper is that doing so, by itself, leads to an approximation error of up to $2$ -- an exact optimum over the original support is only a (possibly tight) $2$-approximation. An exact barycenter for grid-structured data lies in an $N$-times finer grid.  An additional challenge lies in scenarios where measures with sparse support are spread out over a large region, or where they lie in high dimension, such that it is not feasible to discretize the whole underlying space. These restrictions have led to interest in (different) approaches where the support for a barycenter approximation \citep{ccs-18,fms-19,lspc-19} or a discretization of the input measures \citep{scsj-17} is not part of the input. Second, regularization leads to fully dense solutions, which are considered an undesirable `blur' in many applications. This contrasts with the search for sparse exact or approximate barycenters as done in this paper. Post-processing could be used to `sparsify' a dense solution, but it is open whether this can be done in a way to retain a provable approximation guarantee or to obtain a non-mass splitting transport. And third, they sometimes exhibit poor numerical behavior as regularization decreases \citep{kddgtu-19}; typically a fixed number of iterations is hard-coded.

There are other successful algorithms that do not rely on smoothing or regularization: for example, a non-smooth optimization algorithm based on quasi-Newton steps and the fast computation of super-gradients performs well in practice \citep{coo-15}.  Many further examples are based on other types of simplifications of the Wasserstein distance. The so-called Radon barycenters and Sliced barycenters \citep{brpp-15,rpdb-12} are restricted to special instances  (Radon barycenters are only practical for data on a grid, Sliced barycenters deal with support points of uniform mass), but provide good results in low dimension. The idea is to use a Radon transform to obtain $1$-dimensional projections of the support points to lines sampled randomly, from which an expectation of the Wasserstein distance can be devised. Further, a use of the simpler $W_1$-distance instead of the $W_2$-distance leads to the so-called Beckmann problem, which allows various efficient approaches \citep{abgv-19,es-17,srgb-14}.



\subsection{Contributions}\label{sec:contr}

In this paper, we study LP-based approaches to the discrete barycenter problem. In Section \ref{stateoftheart}, we introduce some notation and recall previous related work on linear programming for the problem.  In Section \ref{results}, we present and discuss our main contributions.

First, we show that an optimal measure for Eq.~(\ref{barycenter}), when restricted to the union of supports of the original measures, gives a tight $2$-approximation for the barycenter problem (Theorem \ref{thm:originalsupport}). This result has an immediate implication for algorithms in the literature that compute an approximate barycenter for grid-structured data: if the computation is done over the original grid itself, the algorithm does not converge to an exact barycenter, but to a best approximation of it over the grid, which can give up to a $2$-error. 

Next, we exhibit that a restriction to the support of the original measures allows us to trade this small, provable approximation error for a dramatic improvement in the size of barycenter LPs: we obtain an LP-based $2$-approximation algorithm that can be set up and solved in strongly-polynomial time (Algorithm \ref{algo:2approx}, Theorem \ref{cor:originalsupport}), i.e., polynomial in the number of variables and constraints in the input (and their actual size does not matter). The algorithm finds a {\em sparse} approximation; to the best of the author's knowledge, this is the first algorithm with an approximation guarantee for a sparse solution to the problem. The result shows that the barycenter problem can be efficiently approximated {\em for any data}.


The output of Algorithm \ref{algo:2approx} may not allow for a non-mass splitting transport. (Recall that the existence of such a transport is a property of all exact barycenters and is important for many applications.) Next, we present a second algorithm that improves an approximate barycenter as computed through Algorithm \ref{algo:2approx} to another measure with a non-mass splitting transport of lower cost, and prove that this computation also runs in strongly-polynomial time (Algorithm \ref{algo:nonmasssplit}, Theorems \ref{thm:fixingtheproperties} and \ref{cor:nonmasssplit}, proofs in Appendix \ref{sec:localimprovement}). This algorithm achieves the improvement by moving mass out of the union of original supports to a new, updated support set. Both Algorithms \ref{algo:2approx} and \ref{algo:nonmasssplit} work for any input and do not require the a priori specification of a fixed support set.

Finally, we use the two algorithms as the building blocks of an iterative scheme alternating between them (Algorithm \ref{algo:heuristic}). We prove that it terminates with a $2$-approximation with both sparse support and an associated non-mass splitting optimal transport at the same time (Theorem \ref{thm:heuristic}, proof in Appendix \ref{sec:iterativeproofs}). The theoretical running time of this third algorithm remains open at this time; in practical computations we observe a low number of iterations (often just two to four) before termination. This behavior is reminiscent of the well-known $k$-means algorithm \citep{l-82,m-67}. Further, while we exhibit an example that shows the $2$-approximation bound for Algorithm \ref{algo:heuristic} is tight in theory, we have not observed more than a $20\%$ error, respectively a multiplicative $1.2$ error, in the computations in this paper.

We conclude with sample computations and an analysis of the scaling of the algorithms in Section \ref{computations}. We provide comparisons to exact, LP-based computations, and observe dramatic improvements in running time. We also provide a brief comparison to one of the most popular regularization-based algorithms in the literature (\cite{cd-14}) for grid-structured data. However, one has to be careful in such a comparison; our algorithms have some properties that come at a significant computational cost: Algorithm \ref{algo:2approx} already is an {\em exact, sparse} solution over the original support, and its result is further refined through Algorithms \ref{algo:nonmasssplit} and \ref{algo:heuristic}. Further, the algorithms are numerically stable, work for any data, and work without specification of a fixed support set for the solution. As expected, simplifications to the objective function (like entropy regularized input), and the search for dense approximations over the original support, leads to worse approximation errors, but much faster running times. In fact, we discuss why grid data is an especially poor setting for our methods in view of computational speed: the size of the original support scales quadratically with the density of the grid. 

We close the discussion with a `best-case' example for computational speed: a large number of measures of small, overlapping support (without grid-structure). Due to a linear scaling of the LPs in the number of measures, we are then able to find solutions for thousands of measures. Such data is common in operations research applications that involve geographical locations, but algorithms in the literature are not designed to work (well) in this setting.

\section{Preliminaries}\label{stateoftheart}

We begin by recalling some background on LP-based approaches to the discrete barycenter problem, following \cite{abm-16} and \cite{m-16}. We are given a set of {\em discrete probability measures} $P_1,\dots,P_N$, i.e., they have finite support in $\R^d$ and their total mass sums up to $1$. For a simple wording, we call them {\em measures} in this paper, or {\em discrete measures} to stress that they have a finite support. A set of support points with associated total mass less than $1$ will be called a {\em partial measure}. We denote the {\em support} of $P_i$ as $\supp(P_i)$ and the corresponding number of support points as $|P_i|=|\supp(P_i)|$. $|P_i|$ is called the {\em size} of $P_i$. Further, we are given a weight vector $\lambda=(\lambda_1,\dots,\lambda_N) \in \R^N_{> 0}$ with $\sum_{i=1}^N \lambda_i=1$.

The general definition of a Wasserstein barycenter refers to a measure $\bar P$ on $\Bbb R^d$ which satisfies Eq.~(\ref{barycenter}), i.e.,
\begin{equation*}
\phi(\bar P)=\sum_{i=1}^N \lambda_i W_2(\bar{P},P_i)^2=\inf_{P\in\mathcal P^2(\Bbb R^d)}\sum_{i=1}^N \lambda_i W_2( P, P_i)^2.
\end{equation*}
For discrete measures $P_1,\dots,P_N$, one can show \citep{abm-16} that all optimizers of Eq.~(\ref{barycenter}) must be supported in the finite set $S\subset \mathbb{R}^d$ defined as
\begin{equation}
\label{centers}
S:=\left\{\sum_{i=1}^N \lambda_i x_i : \spc x_i\in\supp(P_i)\right\}.
\end{equation}
$S$ is the set of weighted centroids for all possible combinations of support points, one from each measure $P_i$. Note that $S$ does not have to overlap with the support sets $\supp(P_i)$.

\subsection{Linear programs for discrete barycenters}\label{stateoftheart1}

Setting $\mathcal{P}_{\hspace*{-0.05cm}S}^2(\Bbb R^d):=\{P\in\mathcal{P}^2(\Bbb R^d)|\spc \supp(P)\subseteq S\}$, the infinite-dimensional problem in Eq.~(\ref{barycenter}) can be solved by replacing the requirement $P\in\mathcal P^2(\Bbb R^d)$ with $P\in\mathcal{P}_{\hspace*{-0.05cm}S}^2(\Bbb R^d)$ to obtain
\begin{equation}
\phi(\bar P)=\inf_{P\in\mathcal{P}_{\hspace*{-0.05cm}\mathcal{S}}^2(\Bbb R^d)}\sum_{i=1}^N \lambda_i W_2( P, P_i)^2.
\end{equation}

This yields a finite-dimensional minimization problem, which can be solved through linear programming \citep{abm-16,bs-18,coo-15}. We recall this construction in two steps: we begin with the computation of the value $\phi(P_0)=\sum_{i=1}^N\lambda_i W_2(P_0,P_i)^2$, i.e., the cost of an optimal transport from $P_0$ to all the $P_i$ for a {\em given} $P_0$. Then we make $P_0$ part of the optimization, too.

Let $P_1,\dots,P_N$ be a set of discrete measures and let $\supp(P_i) = \{x_{ik}\big| k = 1,...,|P_i|\}$. Further, let $P_0$ be another (fixed) discrete measure and let $\supp(P_0) = \{x_{j}\big| j = 1,...,|P_0|\}$. Finally, let $d_{ik}$ be the mass of support point $x_{ik}$ in $P_i$ and $d_j$ be the mass of support point $x_j$ in $P_0$. Then we can find the value of $\phi(P_0)$ by solving the following LP: 
\begin{align*}\label{measureLP}\tag{opt. transport}
\min_y   &\spc\spc \sum_{i=1}^N\lambda_i\sum_{j=1}^{|P_0|}\sum_{k = 1}^{|P_i|}\|x_j - x_{ik} \|^2 y_{ijk} \nonumber \\
\text{s.t.} \ \sum_{k=1}^{|P_i|} y_{ijk}  =  &\spc\spc d_j  \spc\spc\spc\spc\forall i=1,\ldots,N,\spc\forall j=1,\ldots,|P_0|\nonumber\\
\sum_{j=1}^{|P_0|} y_{ijk}  =  &\spc\spc d_{ik} \spc\spc\spc\forall i=1,\ldots,N,\spc\forall k=1,\ldots,|P_i|\\
y_{ijk}  \geq  & \spc\spc0 \,\spc\spc\spc\spc\spc\forall i=1,\ldots,N,\spc\forall j=1,\ldots,|P_0|,\spc\forall k=1,\ldots,|P_i|\nonumber
\end{align*}
Note that we not only find the optimal objective function value $\phi(P_0)$, but also a corresponding {\em (optimal) transport} $y=(y_{ijk})_{i\leq N,j \leq |P_0|,k\leq |P_i|}$ between $P_0$ and the $P_1,\dots,P_N$. 

Next, the mass becomes part of the optimization. Instead of just searching for an optimal transport from a fixed measure $P_0$, we use a set $S_0$ of {\em possible support points} with associated variables that represent mass on them. By introducing variables  $z=(z_j)_{j\leq |S_0|}$ for the points in a given set $S_0=\{x_{j}\big| j = 1,...,|S_0|\}$ to denote the possible mass at $x_j\in S_0$, we obtain an LP that both finds an optimal measure $P_0$ supported on $S_0$, as well as a corresponding optimal transport: 


\begin{align}\label{baryLP}\tag{barycenter}
 \min_{z,y}  & \spc\spc \sum_{i=1}^N\lambda_i\sum_{j=1}^{|S_0|}\sum_{k = 1}^{|P_i|}\|x_j - x_{ik} \|^2 y_{ijk}  \nonumber \\
\text{s.t.} \  \sum_{k=1}^{|P_i|} y_{ijk}  = &\spc\spc z_j \spc\spc\spc\spc\forall i=1,\ldots,N,\spc\forall j=1,\ldots,|S_0|\nonumber\\
  \sum_{j=1}^{|S_0|} y_{ijk}  = &\spc\spc d_{ik}  \spc\spc\spc\forall i=1,\ldots,N,\spc\forall k=1,\ldots,|P_i|\nonumber\\
  y_{ijk}  \geq  &\spc\spc0  \,\spc\spc\spc\spc\spc\forall i=1,\ldots,N,\spc\forall j=1,\ldots,|S_0|,\spc\forall k=1,\ldots,|P_i|\nonumber\\
 z_j  \in & \spc\spc\mathbb{R} \spc\spc\spc\spc\spc\forall  j=1,\ldots,|S_0|\spc\nonumber
\end{align}
Note that the variables $z_j$ satisfy $z_j\geq 0$ and $\sum_{j=1}^{|S_0|} z_j = 1$ because of satisfaction of the other constraints and because $\sum_{i=1}^{|P_i|} d_{ik}=1$ for all $i\leq N$. Thus, it suffices to specify $z_j \in \mathbb{R}$.

The above LP computes a measure represented by $z$ and a corresponding optimal transport $y$. For $S_0=S$, the returned $(z,y)$ represents a discrete barycenter by $z$ and a corresponding optimal transport by $y$. For $S_0\neq S$, we call the measure represented by $z$ an {\em $S_0$-barycenter}, an {\em approximation of the barycenter in $S_0$}, or simply an {\em approximate barycenter} when the context is clear.

\subsection{Scaling of the LPs}\label{stateoftheart2}

Let us consider the size of LP (\ref{baryLP}). It consists of $|S_0|+ |S_0|\cdot \sum_{i=1}^N |P_i|$ variables and $N\cdot |S_0| + \sum_{i=1}^N |P_i|$ equality constraints. For the computation of an exact barycenter, we set $S_0=S$. In this case, we get a worst-case bound of $|S_0|=\prod_{i=1}^N |P_i|$. Let now $|P_{\max}|=\max_{i=1,\dots,N} |P_i|$. If all measures have the same number of support points, we get $\sum_{i=1}^N |P_i| = N\cdot |P_{\max}|$ and $\prod_{i=1}^N |P_i|=|P_{\max}|^N$. So we have an LP of up to $|P_{\max}|^N +|P_{\max}|^N \cdot N\cdot|P_{\max}|$ variables and  $N\cdot |P_{\max}|^N + N \cdot |P_{\max}|$ equality constraints.

A refined analysis reveals that some of the variables and constraints can be redundant. For example, if the measures overlap in some of their support points, then $|S_0|$ and consequently the size of the LP becomes smaller. In fact, LP (\ref{baryLP}) is always of polynomial size for data on a grid \citep{bs-18}. However, in general one cannot rule out a scaling of the size of the LP for $S_0=S$ that is exponential in $N$ even if $|P_{\max}|$ is fixed, and a polynomial scaling in $|P_{\max}|$ even if $N$ is fixed. The main reason why it was possible to compute an exact barycenter for the example in \cite{abm-16} with only $8$ measures of $9$ support points was the fact that all measures had the same small support, which had a dramatic effect in reducing $|S_0|$. This highlights the potential benefit from performing an approximate computation where one replaces $S$ by a smaller set $S_0$.

\subsection{Sparsity and non-mass split}\label{stateoftheart3}

The feasible regions of LPs (\ref{measureLP}) and (\ref{baryLP}) are bounded, and thus standard arguments of linear programming show that there always exists an optimal vertex. In a vertex, an inclusion-maximal set of variables is set to $0$. By a careful analysis of which of the variables $z_j,y_{ijk} $ are equal to $0$, it is possible to show a first favorable property: in contrast to the large number $|S|$ of possible support points, which can be up to $\prod_{i=1}^N |P_i|$, there always exists a barycenter that assigns nonzero mass to at most $\sum_{i=1}^N |P_i|  - N + 1$ points \citep{abm-16}.

\begin{proposition}\label{sparsethm} Let  $P_1,\ldots,P_N$ be discrete measures. Then for any weights $\lambda\in\R^n_{>0}$, there exists a barycenter $\bar{P}$ of these measures such that the size $|\bar{P}|$ satisfies
\begin{equation}\label{sparseequ}
|\bar{P}| \leq \sum_{i=1}^N |P_i|  - N + 1.
\end{equation}
\end{proposition}

We call a measure $\bar P$ that satisfies $|\bar{P}| \leq \sum_{i=1}^N |P_i|  - N + 1$ {\em sparse}. Proposition \ref{sparsethm} states that there always exists a {\em sparse barycenter}. A proof is based on the existence of an optimal vertex of the polyhedron for LP (\ref{baryLP}) \citep{abm-16,m-16}. The argument also works if a support set $S_0\neq S$ is used. LP (\ref{baryLP}) then optimizes the objective function in Eq.~(\ref{barycenter}) over the set $\mathcal{P}_{S_0}^2(\R^d)$  of all measures $P$ with support in $S_0$. For these different support sets, we have the following generalization of Proposition \ref{sparsethm}.

\begin{corollary}\label{sparsecor} Let  $P_1,\ldots,P_N$ be discrete measures in $\mathbb{R}^d$, let $S_0=\{x_j: j=1,\dots,|S_0|\} \subset \mathbb{R}^d$, and let  $\mathcal{P}_{S_0}^2(\R^d)$  be the set of all measures $P$ with support in $S_0$. Then for any weights $\lambda\in\R^n_{>0}$, there exists an approximate barycenter $\bar{P}_0$ in $S_0$ such that the size $|\bar{P}_0|$ satisfies
\begin{equation}\label{sparseequ2}
|\bar{P}_0| \leq \sum_{i=1}^N |P_i|  - N + 1.
\end{equation}
\end{corollary}


Further, for any exact barycenter $\bar P$ there exists a {\em non-mass splitting optimal transport} from $\bar P$ to the $P_1,\dots,P_N$ \citep{abm-16,m-16}. This means that for all $x_j\in \supp(\bar P)$ with mass $d_j$ and for each $i$, there is exactly one $k$ with $y_{ijk}=d_j$, while $y_{ijk'}=0$ for all $k'\neq k$. Each support point of a barycenter only transports mass to exactly one support point in each measure. In this case, we say that a support point {\em does not split mass} or that a support point is {\em non-mass splitting}. 

In fact, {\em any } optimal transport from a discrete barycenter $\bar P$ to the corresponding set of measures is non-mass splitting. While this has not been stated explicitly in \citep{abm-16}, it is not hard to prove: recall that the (weighted) centroid $c$ of a set of points $x_1,\dots,x_n$ is the unique minimizer of a functional that measures the (weighted) summed-up squared Euclidean distances of a single point to all points in the set. This can be seen through a simple transformation \begin{eqnarray*}\sum\limits_{i=1}^N \lambda_i \|(s+c)-x_i\|^2= s^Ts-c^Tc+\sum\limits_{i=1}^N \lambda_i x_i ^Tx_i,
\end{eqnarray*}
which is minimal for $s^Ts=0$, so $s=0$. If there was a barycenter support point splitting mass, it could be split into two (or more) centroids of support points in the measures of the same total mass, and the cost of transport would be strictly lower. We formally state this observation. 

\begin{proposition}\label{centroid} Let $P_1,\ldots,P_N$ be discrete measures, and let $\bar P$ be a barycenter for these measures. Then any optimal transport from $\bar P$ to $P_1,\ldots,P_N$ is non-mass splitting.
\end{proposition}

\section{Main Results}\label{results} 

In this paper, we study approximations of the discrete barycenter problem where the set $S$, required to find an exact barycenter, is replaced by a much smaller set $S_0$. This is motivated by the unfavorable scaling of LP (\ref{baryLP}) with respect to $|S|$, respectively $|S_0|$; see Section \ref{stateoftheart2}. 

\subsection{A strongly-polynomial $2$-approximation}\label{sub2-approx}

Recall that the set of possible support points of a discrete barycenter is
\begin{equation}
S:=\left\{\sum_{i=1}^N \lambda_i x_i : \spc x_i\in\supp(P_i)\right\},
\end{equation}
which may consist of up to $\prod\limits_{i=1}^N |P_i|$ points. This is a much larger number than the size of the union of supports of the measures
\begin{equation}\label{eq:S_org}
S_{\text{org}}:=\bigcup\limits_{i=1}^N \supp(P_i),
\end{equation}
which satisfies $|S_{\text{org}}|\leq \sum\limits_{i=1}^N |P_i|$ with equality if and only if the supports are disjoint. 

Note that the maximal size of $S_{\text{org}}$ only barely exceeds the bound given in Proposition \ref{sparsethm}. First, we show that the approximation error from searching for an approximate barycenter in $S_{\text{org}}$, i.e., setting $S_0=S_{\text{org}}$ in LP (\ref{baryLP}), can be bounded by a factor of two. This bound is tight.


\begin{theorem}\label{thm:originalsupport}
Let $\bar P$ be a barycenter and let $\bar P_\text{org}$ be an approximate barycenter in $S_{\text{org}}$.
Then 
\begin{equation*}
\phi(\bar P_\text{org}) \leq 2 \cdot \phi(\bar P) 
\end{equation*}
and this bound can become tight, i.e., there is a set of measures $P_1,\dots,P_N$ and a set of weights $\lambda_1,\dots,\lambda_N$ for which $\phi(\bar P_\text{org}) = 2 \cdot \phi(\bar P)$.
\end{theorem}

\begin{proof}
We denote the mass of a support point $c$ of a barycenter $\bar P$ by $d_c$. By Proposition \ref{centroid}, there is an optimal transport such that $c$ transports its mass to exactly one support point $x_i$ in each $P_i$ for all $i\leq N$. Due to optimality of $\bar P$, $c$ is the weighted centroid $c=\sum_{i=1}^N \lambda_i x_i$ of these points. Recall the discussion after Corollary \ref{sparsecor}.

Each support point $c$ contributes $d_c\cdot \sum_{i=1}^N \lambda_i \|c-x_i\|^2$ to the corresponding value $\phi(\bar P)$. Let $s\in S_{\text{org}}=\bigcup_{i=1}^N \supp(P_i)$ be such that $\|s-c\|^2$ is minimal and note that
\begin{eqnarray*}\sum\limits_{i=1}^N \lambda_i \|s-x_i\|^2 = s^Ts -2c^Ts + \sum\limits_{i=1}^N\lambda_i x_i^Tx_i =
 (s^Ts -2c^Ts + c^Tc) + \\ + (c^Tc -2c^Tc+  \sum\limits_{i=1}^N\lambda_i x_i^Tx_i) =  
 \sum\limits_{i=1}^N\lambda_i (\|s-c\|^2+\|c-x_i\|^2)\end{eqnarray*}
for any $s$. By choice of $s$ and the fact that $x_i\in \supp(P_i)$, we know $\|s-c\|^2\leq \|c-x_i\|^2$ for all $i\leq N$, so we obtain
$$\sum\limits_{i=1}^N \lambda_i \|s-x_i\|^2=\sum\limits_{i=1}^N \lambda_i (\|s-c\|^2+\|c-x_i\|^2)\leq 2 \cdot \sum\limits_{i=1}^N \lambda_i \|c-x_i\|^2.$$
Thus the transport from $s$, instead of from $c$ itself, introduces an approximation error of $2$, i.e., each such $s$ contributes at most $2\cdot d_c  \sum_{i=1}^N \lambda_i \|c-x_i\|^2$ to the value $\phi(\bar P_\text{org})$. As this holds for all weighted centroids $c\in \supp(\bar P)$ and corresponding closest $s \in S_{\text{org}}$, this shows the existence of a measure $\bar P_\text{org} \in \mathcal{P}_{\text{org}}^2(\R^d)$ with approximation error $2$ with respect to $\phi$.

It remains to prove that the bound can be tight. We do so through a simple example. Let $P_1, P_2$ be two measures with a single support point $x_{11}\in \supp(P_1)$, $x_{21}\in \supp(P_2)$, each of mass $1$. Then $\bar P$ consists of the single support point $c=\lambda_1 x_{11} + \lambda_2 x_{21}$ of mass $1$ and thus 
\begin{eqnarray*}\phi(\bar P)=\lambda_1 \cdot \|c- x_{11}\|^2 + \lambda_2 \cdot \|c- x_{21}\|^2 = \lambda_1\cdot \|\lambda_2 (x_{21}-x_{11})\|^2 +  \lambda_2\cdot \|\lambda_1 (x_{11}-x_{21})\|^2= \\ = \lambda_1\lambda_2(\lambda_2+\lambda_1)\|x_{21}-x_{11}\|^2=\lambda_1\lambda_2\|x_{21}-x_{11}\|^2.\end{eqnarray*}
In contrast, the restriction of an approximate barycenter $\bar P_\text{org}$ to possible support $S_{\text{org}}=\{x_{11},x_{21}\}$ would give $\phi(\bar P_\text{org})=\min\{\lambda_1, \lambda_2\} \cdot \|x_{21}-x_{11}\|^2$. Note $\lambda_1\cdot \lambda_2\geq \frac{1}{2} \min\{\lambda_1, \lambda_2\}$, with equality if and only if $\lambda_1 = \lambda_2=\frac{1}{2}$. In this case, $\phi(\bar P_\text{org})=2\cdot \phi(\bar P)$. \hfill$\square$ 
\end{proof}

\begin{figure}
\hspace*{-0.5cm}
\subfloat[Measure $P_1$]{
  \includegraphics[scale=0.28]{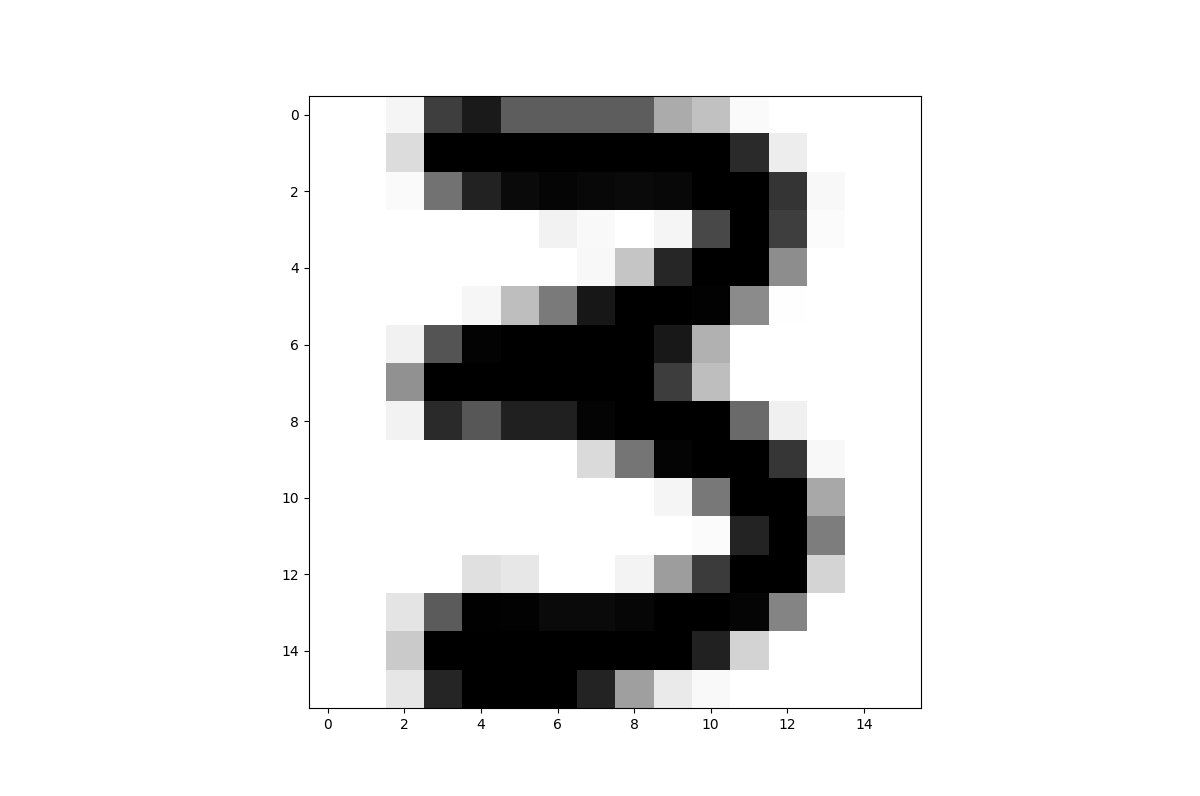}
}\hspace*{-2cm}
\subfloat[Measure $P_2$]{
  \includegraphics[scale=0.28]{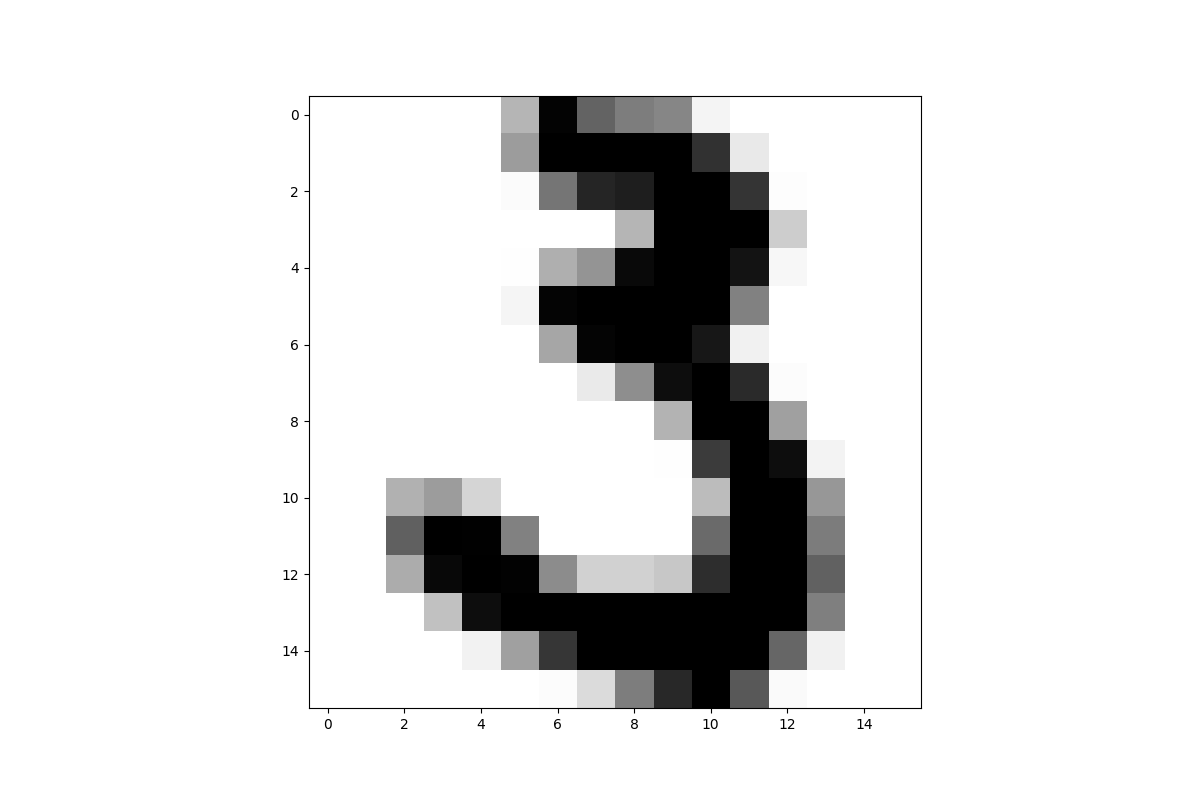}
}

\hspace*{-0.5cm}
\subfloat[Measure $P_3$]{
  \includegraphics[scale=0.28]{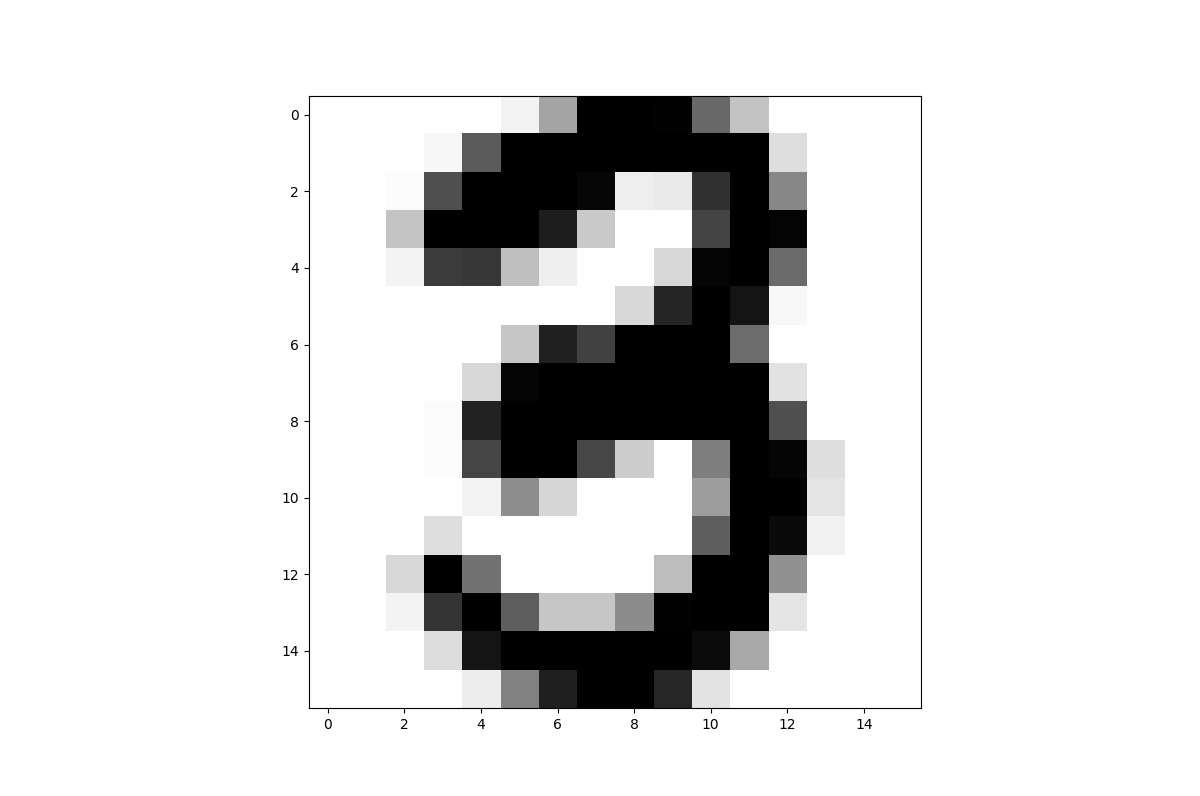}
}\hspace*{-2cm}
\subfloat[Measure $P_4$]{
  \includegraphics[scale=0.28]{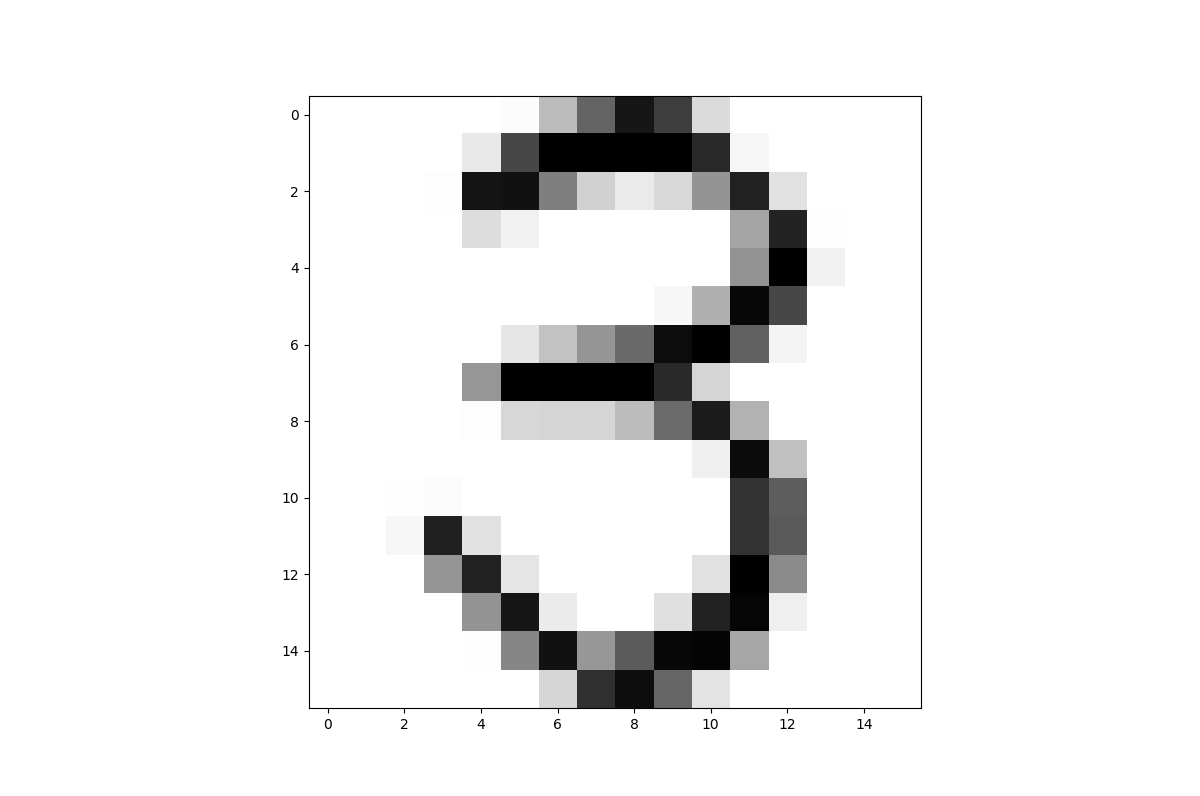}
}

\hspace*{-0.5cm}
\subfloat[Barycenter $\bar P$]{
  \includegraphics[scale=0.28]{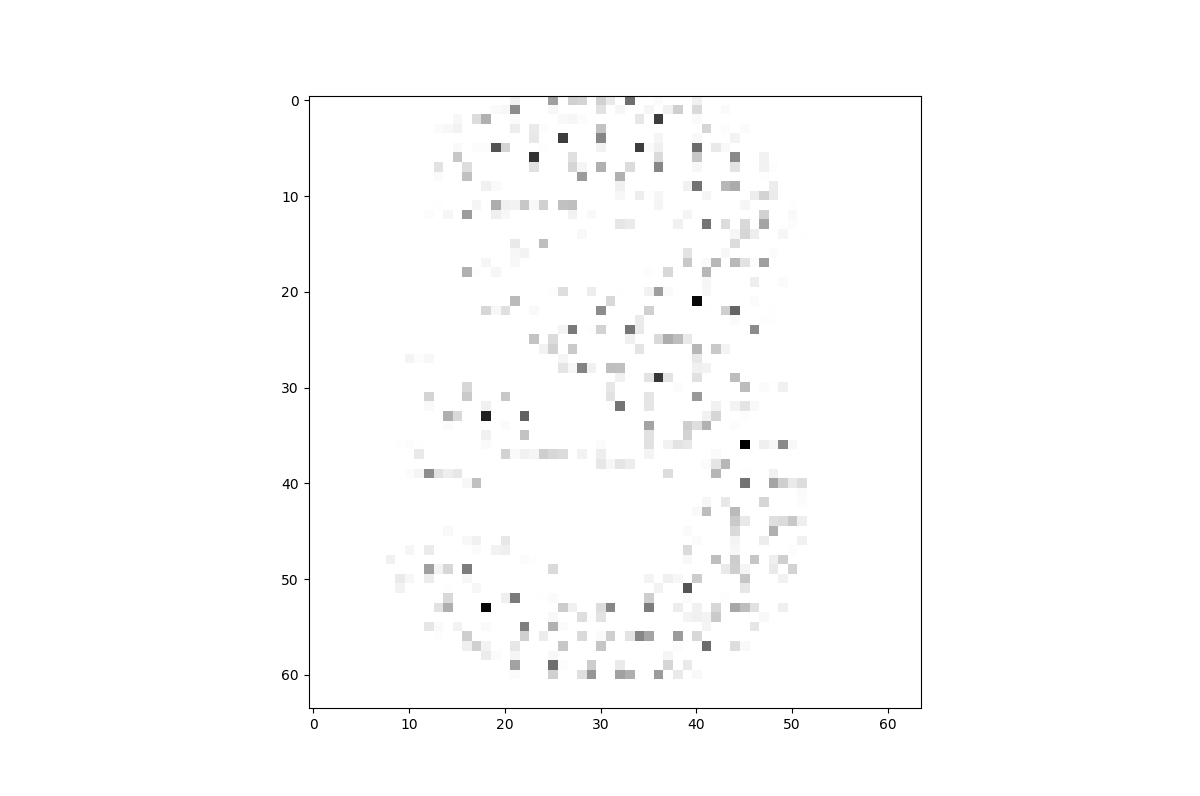}
}\hspace*{-2cm}
\subfloat[Approximate Barycenter $\bar P_\text{org}$]{
  \includegraphics[scale=0.28]{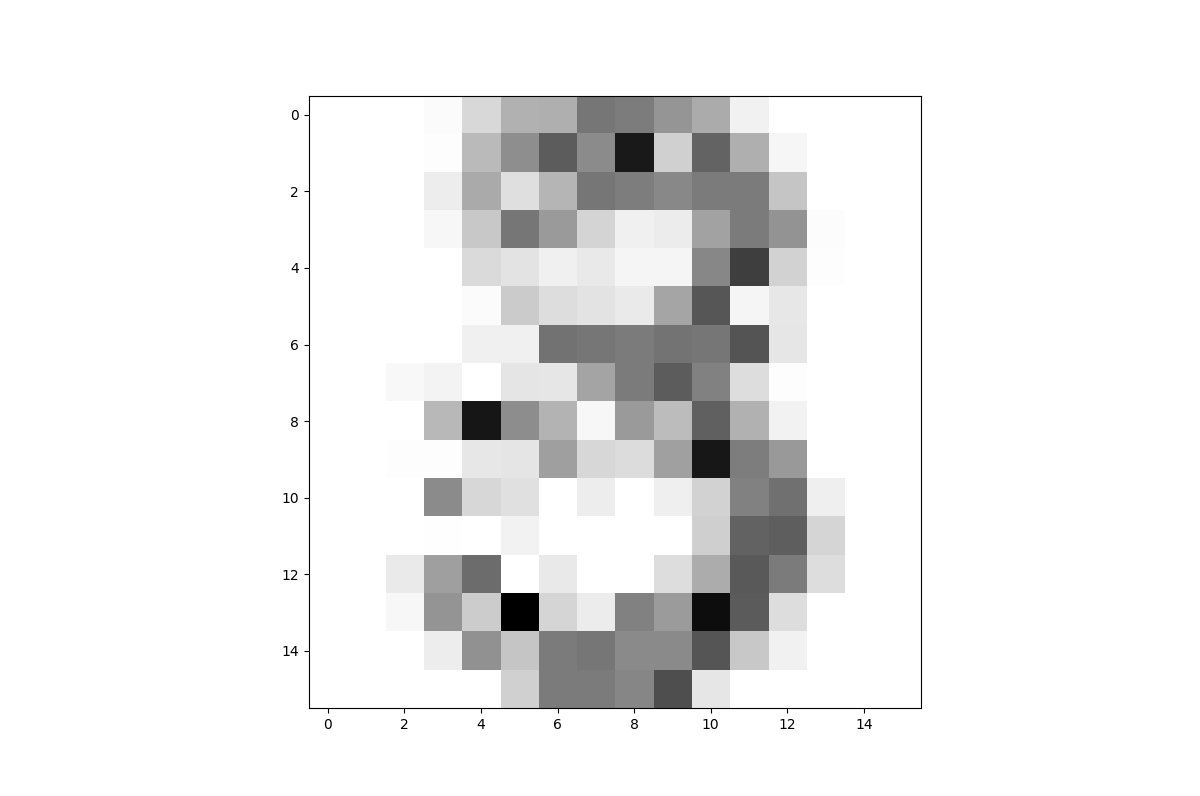}
}
\caption{Four measures $P_1,\dots,P_4$ supported on a $16 \times 16$ grid in the first two rows. The bottom row shows a barycenter $\bar P$ and an approximate barycenter $\bar P_\text{org}$. While the support of $\bar P_\text{org}$ lies in the original $16 \times 16$ grid, the support for $\bar P$ lies in a four times finer grid.} 
\label{fig:differentsupports}
\end{figure}

The difference between the support for an exact barycenter and for an approximation in $S_\text{org}$ is highlighted in Figure \ref{fig:differentsupports}: the first two rows show four handwritten digits scanned into a $16 \times 16$ grid. (See \cite{lbh-98} for some information on this data set.) These are the measures $P_1,\dots,P_4$. The varying shades of grey indicate different masses at the support points of the grid (the darker, the larger the mass). The masses for each measure add up to $1$. The bottom row depicts an exact barycenter and a $2$-approximation in the original $16 \times 16$ grid (for all $\lambda_i=\frac{1}{4}$). The support grid for the exact barycenter is four times finer, a $(4\cdot 16-3)\times(4\cdot 16-3)=61\times 61$ grid.

\begin{algorithm}
\setstretch{1.0}
\vspace{0.2cm}

 {\bf Input}
\begin{itemize}
\item Measures $P_1,\dots,P_N \subset \mathbb{R}^d$, support $S_\text{org}=\bigcup_{i=1}^N \supp(P_i)$
\item $\lambda_1,\dots,\lambda_N> 0$ with $\sum_{i=1}^N \lambda_i =1$
\end{itemize}
 {\bf Algorithm}
\\\\ \noindent Compute an approximate barycenter $\bar P_\text{org}$ in $S_{\text{org}}$ as an optimal vertex $(z,y)$ of
\begin{align} 
\min_{z,y}  &\spc\spc\phi(\bar P_\text{org}) := \sum_{i=1}^N\lambda_i\sum_{j=1}^{|S_{\text{org}}|}\sum_{k = 1}^{|P_i|}\|x_j - x_{ik} \|^2 y_{ijk} \nonumber \\
\text{s.t.} \ \sum_{k=1}^{|P_i|} y_{ijk}  = &\spc\spc z_j  \spc\spc\spc\spc\forall i=1,\ldots,N,\spc\forall j=1,\ldots,|S_{\text{org}}|\nonumber\\
\sum_{j=1}^{|S_{\text{org}}|} y_{ijk}  = &\spc\spc d_{ik}\spc\spc\spc\forall i=1,\ldots,N,\spc\forall k=1,\ldots,|P_i|\nonumber\\
y_{ijk}  \geq & \spc\spc0 \,\spc\spc\spc\spc\spc\forall i=1,\ldots,N,\spc\forall j=1,\ldots,|S_{\text{org}}|,\spc\forall k=1,\ldots,|P_i|\nonumber\\
 z_j  \in & \spc\spc\mathbb{R} \spc\spc\spc\spc\spc\forall  j=1,\ldots,|S_{\text{org}}|\spc\nonumber
\end{align}
and return $z$ to represent $\bar P_\text{org}$ and the corresponding optimal transport $y$.
  \caption{Sparse $2$-approximate barycenter in the original support}
   \label{algo:2approx}
\end{algorithm}

We formally denote the choice of $S_\text{org}$ in LP (\ref{baryLP}), as performed for Theorem \ref{thm:originalsupport}, as Algorithm \ref{algo:2approx}. Note that the algorithm is stated to compute an optimal {\em vertex} of the feasible region. (For a convenient wording, we will say that Algorithm \ref{algo:2approx} is used with a given different support $S_0$ as input when we require an optimal vertex, and not just any optimal solution, of LP (\ref{baryLP}).) The search for a vertex guarantees that the sparsity condition stated in Corollary \ref{sparsecor} is satisfied, so the returned measure is not only an approximate barycenter in $S_\text{org}$, but also {\em sparse}. However, it is possible that any corresponding optimal transport splits mass (which would not happen for an exact barycenter).  Here is an example of this type. 


\begin{example}\label{ex1}
Consider the two measures $P_1, P_2$ depicted at the top of Figure \ref{fig:masssplit}. The different radii of the filled circles represent the different masses on the support points.  Let $\lambda_1=\lambda_2=\frac{1}{2}$

Measure $\bar P_\text{org}\in \mathcal{P}_{\text{org}}^2(\R^d)$ (second row, left) is an optimal barycenter approximation in $S_{\text{org}}$. It consists of only two support points, while $P_1$ and $P_2$ have three support points. Thus, there exists a support point of $\bar P_\text{org}$ that splits mass in any transport, including the unique optimal one (second row, right): the top support point of $\bar P_\text{org}$ transports its mass $\frac{1}{2}$ in two parts $\frac{1}{4}$ to two support points of $P_1$; the same happens in the bottom part with respect to $P_2$. Such a split of mass does not happen for an exact barycenter (third row). \hfill$\square$

\begin{figure}[t]
\begin{center}
\subfloat[Measure $P_1$]{
\begin{tikzpicture}[scale=1.5]
\node (11) at (0,1.5) {$\frac{1}{4}$};
\node (12) at (1,-0.5) {$\frac{1}{2}$};
\node (21) at (2,1.5) {$\frac{1}{4}$};

\fill [blue] (0,1) circle (3pt);
\fill [blue] (1,0) circle (4.5pt);
\fill [blue] (2,1) circle (3pt);
\fill [black] (0,0) circle (1pt);
\fill [black] (1,1) circle (1pt);
\fill [black] (2,0) circle (1pt);
\end{tikzpicture}
}
\qquad \qquad  \qquad  \qquad  \qquad  \qquad \;
\subfloat[Measure $P_2$]{
\begin{tikzpicture}[scale=1.5]
\node (11) at (0,-0.5) {$\frac{1}{4}$};
\node (12) at (1,1.5) {$\frac{1}{2}$};
\node (21) at (2,-0.5) {$\frac{1}{4}$};

\fill [black] (0,1) circle (1pt);
\fill [black] (1,0) circle (1pt);
\fill [black] (2,1) circle (1pt);
\fill [red] (0,0) circle (3pt);
\fill [red] (1,1) circle (4.5pt);
\fill [red] (2,0) circle (3pt);
\end{tikzpicture}
}

\subfloat[Measure $\bar P_\text{org}$]{\label{fig:test}
\begin{tikzpicture}[scale=1.5]
\node (11) at (1,-0.5) {$\frac{1}{2}$};
\node (12) at (1,1.5) {$\frac{1}{2}$};

\node (dummy) at (1,-0.75) {$ $};

\fill [black] (1,0) circle (4.5pt);
\fill [black] (1,1) circle (4.5pt);

\fill [black] (0,1) circle (1pt);
\fill [black] (1,0) circle (1pt);
\fill [black] (2,1) circle (1pt);
\fill [black] (0,0) circle (1pt);
\fill [black] (1,1) circle (1pt);
\fill [black] (2,0) circle (1pt);
\end{tikzpicture}} 
\qquad \qquad  \qquad  \qquad  \qquad  \qquad  \;\;\;
\subfloat[Mass Split]{
\begin{tikzpicture}[scale=1.5]
\node (11) at (1,-0.7) {$ $};

\fill [blue] (0,1) circle (3pt);

\fill [blue] (2,1) circle (3pt);
\fill [red] (0,0) circle (3pt);

\fill [red] (2,0) circle (3pt);

\draw[black,line width=0.6mm,->] (1,0) -- (0.05,0);
\draw[black,line width=0.6mm,->] (1,0) -- (1.95,0);

\draw[black,line width=0.6mm,->] (1,1) -- (0.05,1);
\draw[black,line width=0.6mm,->] (1,1) -- (1.95,1);

\fill [blue] (1,0) circle (4.5pt);
\fill [red] (1,1) circle (4.5pt);
 \path[black,line width=0.7mm,->] (1,1.1)   edge[out=135,in=45, my loop] node  {} (1,1.1);
 \path[black,line width=0.7mm,->] (1,-0.1)   edge[out=235,in=315, my loop] node  {} (1,-0.1);

\end{tikzpicture}}

\subfloat[Barycenter $\bar P$]{
\begin{tikzpicture}[scale=1.5]
\node (33) at (1,-0.3) {$ $};

\node (44) at (1,0.7) {$ $};

\node (11) at (0.25,.75) {$\frac{1}{4}$};
\node (12) at (1.25,0.75) {$\frac{1}{2}$};
\node (21) at (2.25,0.75) {$\frac{1}{4}$};

\fill [black] (0,0.5) circle (3pt);
\fill [black] (1,0.5) circle (4.5pt);
\fill [black] (2,0.5) circle (3pt);

\fill [black] (0,1) circle (1pt);
\fill [black] (1,0) circle (1pt);
\fill [black] (2,1) circle (1pt);
\fill [black] (0,0) circle (1pt);
\fill [black] (1,1) circle (1pt);
\fill [black] (2,0) circle (1pt);
\end{tikzpicture}}
\qquad \qquad  \qquad  \qquad  \qquad  \qquad
\subfloat[no Mass Split]{
\begin{tikzpicture}[scale=1.5]
\node (11) at (1,-0.3) {$ $};

\node (22) at (1,0.7) {$ $};

\fill [black] (0,1) circle (1pt);
\fill [black] (1,0) circle (1pt);
\fill [black] (2,1) circle (1pt);
\fill [black] (0,0) circle (1pt);
\fill [black] (1,1) circle (1pt);
\fill [black] (2,0) circle (1pt);

\fill [black] (0,0.5) circle (3pt);
\fill [black] (1,0.5) circle (4.5pt);
\fill [black] (2,0.5) circle (3pt);

\fill [blue] (1,0) circle (4.5pt);
\fill [red] (1,1) circle (4.5pt);

\fill [blue] (0,1) circle (3pt);

\fill [blue] (2,1) circle (3pt);
\fill [red] (0,0) circle (3pt);

\fill [red] (2,0) circle (3pt);

\draw[black,line width=0.6mm,->] (0,0.5) -- (0,0.05);
\draw[black,line width=0.8mm,->] (1,0.5) -- (1,0.95);
\draw[black,line width=0.6mm,->] (2,0.5) -- (2,0.05);

\draw[black,line width=0.6mm,->] (0,0.5) -- (0,0.95);
\draw[black,line width=0.8mm,->] (1,0.5) -- (1,0.05);
\draw[black,line width=0.6mm,->] (2,0.5) -- (2,0.95);
\end{tikzpicture}}

  \end{center}
\caption{Two measures $P_1$, $P_2$ in the top row. An optimal approximate barycenter $\bar P_\text{org} \in \mathcal{P}_{\text{org}}^2(\R^2)$ and the corresponding mass splitting transport in the second row. The exact barycenter and a corresponding non-mass splitting transport in the third row.} 
\label{fig:masssplit}
\end{figure}
\end{example}

We would like to note that the $2$-bound in Theorem \ref{thm:originalsupport} can only be tight in very special cases. Let $s\in S_{\text{org}}$ be such that $\|s-c\|^2$ is minimal for a given weighted centroid $c\notin S_\text{org}$ transporting to $x_{i1},\dots,x_{N1}$ with $x_{i1}\in P_i$ . Then the approximation error $2$ is not tight if $\|c-x_{i1}\|^2 \neq \|c-x_{j1}\|^2$ for any $i\neq j$. This holds because then 
$$\sum\limits_{i=1}^N \lambda_i \|s-x_{i1}\|^2=\sum\limits_{i=1}^N \lambda_i (\|c-x_{i1}\|^2+\|s-c\|^2)< 2 \cdot \sum\limits_{i=1}^N \lambda_i \|c-x_{i1}\|^2,$$
as there has to be an $i\leq N$ with $\|c-x_{i1}\|^2 > \|s-c\|^2$. Further, it is easy to give examples where even an exact barycenter is actually contained in $S_\text{org}$. For example, take arbitrary measures $P_1,\dots,P_N$, compute their barycenter $\bar P$, and consider a new collection $P_1,\dots,P_N,P_{N+1}$ with $P_{N+1}=\bar P$. Then one has $\phi(\bar P_\text{org})=\phi(\bar P)$.



Next, we prove that Algorithm \ref{algo:2approx} runs in strongly-polynomial time. Recall that LPs are generally solvable in weakly-polynomial time, i.e., the number of arithmetic operations is polynomial in the length of a bit representation of the input. This means polynomiality in the number of variables and constraints, as well as in the logarithm of absolute values of numbers in the input. In contrast, a strongly-polynomial running time restricts polynomiality to only the number of variables and constraints.

\begin{theorem}\label{cor:originalsupport}
For all rational input, a $2$-approximate barycenter can be computed in strongly-polynomial time.
\end{theorem}

\proof{}
A proof of strong polynomiality for Algorithm \ref{algo:2approx} is based on exhibiting that LP (\ref{baryLP}) is of strongly-polynomial size, and that its parameters can be computed in strongly-polynomial time. General LPs are known to be solvable in weakly-polynomial time. However, it suffices to restrict the dependency of the running time only to the parameters that appear in the constraint matrix; the numbers in the objective function or the right-hand side of the constraints do not matter \citep{t-86}.



First, note that the constraint matrix of LP (\ref{baryLP}) for $S_{\text{org}}$ only consists of entries in $\{-1,0,1\}$. For the claim of strongly-polynomial solvability, it only remains to prove that the number of variables and constraints of the LP is strongly-polynomial in the size of the input, and that the parameters that appear in the objective function and right-hand sides can be computed from the original input in strongly-polynomial time. 

Let $\mathcal{I}$ be an instance of the problem and let $|\mathcal{I}|$ be the number of bits to represent the input. Any representation of the input  $\mathcal{I}$ has to satisfy  $|\mathcal{I}|\geq \sum_{i=1}^N |P_i|$. As $|S_0|=|S_\text{org}|\leq \sum_{i=1}^N |P_i| \leq |\mathcal{I}|$, LP (\ref{baryLP}) indeed has a strongly-polynomial number of constraints and variables.

The actual numbers that appear in the LP are of types $\lambda_i$, $d_{ik}$, or $\|x_j-x_{ik}\|^2$. The $\lambda_i$ and $d_{ik}$ appear directly in the input, and so do the vectors $x_j$ and $x_{ik}$. As we use rational input, $\|x_j-x_{ik}\|^2=(x_j-x_{ik})^T(x_j-x_{ik})$ is a rational number derived by the sum over products of pairs of coefficients in $x_j$ and $x_{ik}$. This implies that $\|x_j-x_{ik}\|^2$ can be computed in strongly-polynomial time (polynomial in $\text{log }x_j +\text{log }x_{ik}$) and represented in a number of bits that is strongly-polynomial in the number of bits of the original representation of  $x_j, x_{ik}$. This proves the claim.  \hfill$\square$ 
\endproof

\subsection{Recovery of Non-Mass Split}

Next, we design an algorithm that begins with a (sparse) $2$-approximate barycenter computed by Algorithm \ref{algo:2approx}. The algorithm improves it to another measure supported on a subset of $S$ (instead of $S_\text{org}$), for which there exists a {\em non-mass splitting} transport of lower cost, i.e., the approximation error can only become better. Algorithm \ref{algo:nonmasssplit} sums up the approach in pseudocode. We here describe the algorithm in some detail; additional technical details are given in the proof of Theorem \ref{thm:fixingtheproperties} in Appendix \ref{sec:localimprovement}.

 The algorithm greedily breaks up each support point (that splits mass) of the approximate barycenter into several non-mass splitting support points (Steps $1-3$). In the end, all of the non-mass splitting support points are combined to a new measure (Step $4$). The preprocessing performed in Step $2$ guarantees that the non-mass split property for each support point in Step $3$ transfers to a non-mass splitting transport for the new measure constructed in Step $4$. Figure \ref{fig:exampleforstep2} shows a run of the algorithm, which is discussed in more detail as Example \ref{ex0} at the end of the section. 

\noindent{\bf Step 1.} First, the approximate barycenter $\bar P_\text{org}$ and measures $P_1,\dots,P_N$ are broken up into disjoint parts; each part corresponds to a support point $s_l=x_{t_l}$ in the approximate barycenter. By construction, each $P_i^l$ consists of those support points in $P_i$ to which $s_l$ transports mass. The mass of a support point in $P_i^l$ equals the mass it receives as transport from $s_l$. Then we assign new indices to the support points in $P_i^l$ and their masses for a simpler notation, so we do not have to refer to the original $z$ or $y$ in the other steps.

\noindent{\bf Lexicographic Ordering.}  Step $2$ and Step $3$ are based on the construction of so-called {\em lexicographically maximal} vectors. A vector $a=(a_1,\dots,a_n)$ is {\em lexicographically larger} than a vector $b=(b_1,\dots,b_n)$ if there is an index $j\leq n$ such that $a_j > b_j$, and $a_i \geq b_i$ for all $i<j$. For example, the vector $a=(2,2,0,1)$ is lexicographically larger than $b=(2,1,5,10)$. Lexicographic maximality with respect to a set states that there is no lexicographically larger vector in the set. Note that the term gives rise to a total ordering.

The intuition for the construction of lexicographically maximal vectors is to resolve ties. This is necessary in two different settings: in Step 2, as much mass as possible is greedily shifted to support points of lower indices; in Step 3, a lexicographically decreasing sequence of weighted centroids is created from each support point. Together, these two steps make sure that all the weighted centroids that are merged to form $\bar P'$ in Step 4 are distinct and only transport to a single support point in each measure, implying the existence of a non-mass splitting transport.

\noindent{\bf Step 2.} Step $2$ iteratively transforms $(d_1,\dots,d_r)$ to be lexicographically larger and larger while retaining an approximate barycenter supported in $\supp(\bar P_\text{org})$ (that is, the cost of an optimal transport does not increase). It does so via a greedy scheme, where mass is moved to support points in $\supp(\bar P_\text{org})$ with the lowest indices, until this is not possible anymore. We call a $(d_1,\dots,d_r)$ that is not altered by Step $2$ (anymore) {\em greedily lexicographically maximal}. Note that such a vector need not be lexicographically maximal among {\em all} approximate barycenters with the same support, but this is enough for our purposes.

The two loops for $l$ and $j$ establish an order for checking whether mass can be moved from $s_l$ to $s_j$ while keeping optimality over $\supp(\bar P_\text{org})$. The indices $q_i= \text{arg}\max_{q \leq |P_i^l|} (s_j-s_l)^Tx^l_{iq}$ selected in $2a)$ identify support points in the $P_i^l$ that lie the furthest in direction of $s_j-s_l$. 

\begin{algorithm}[H]
\vspace{0.2cm}
\small 
 {\bf Input}
\begin{itemize}
\item Measures $P_1,\dots,P_N \subset \mathbb{R}^d$
\item (sparse) 2-approximate barycenter $\bar P_\text{org}$ and an optimal transport $(z,y)$ (from Alg. \ref{algo:2approx})
\item $\lambda_1,\dots,\lambda_N> 0$ with $\sum_{i=1}^N \lambda_i =1$
\end{itemize}
 {\bf Algorithm}
\begin{enumerate}
\item  {\bf \normalsize (Break up $\bar P_\text{org}$ and $P_1,\dots,P_N$ into parts for each support point of $\bar P_\text{org}$)}\\
\noindent Let $\supp(\bar P_\text{org})=\{s_1,\dots,s_r\}=\{x_{t_1},\dots,x_{t_r}\}$ with corresponding masses $d_1=z_{t_1},\dots,d_r=z_{t_r}$. \\For each $l\leq r$ and $i\leq N$, construct $P_i^l$ (a set of support points with masses) by the rule: 
{\begin{center}
$y_{it_lk}>0$ $\;\;\Rightarrow\;\;$ add $x_{ik}$ to $\supp(P_i^l)$ with mass $y_{it_lk}$
\end{center}}
\noindent Now assign indices for the $P_i^l$ to obtain a notation $P_i^l=\{x^l_{i1},\dots,x^l_{i|P_i^l|}\}$ with corresponding masses $d^l_{i1},\dots,d^l_{i|P_i^l|}$  for all $l\leq r$ and $i\leq N$. 

\item {\bf \normalsize (Make masses $(d_1,\dots,d_r)$ greedily lexicographically maximal)}\\
\noindent For $l=r$ descending to $l=1$ \\
\noindent  \hspace*{0.5cm} For $j=1$ ascending to $j=l-1$\\
\noindent  \hspace*{1cm} $a)$ For each $i\leq N$, identify an index $q_i= \text{arg}\max_{q \leq |P_i^l|} (s_j-s_l)^Tx^l_{iq}$. Then compute\\\noindent  \hspace*{1.4cm}  the weighted centroid $c=\sum_{i=1}^N \lambda_i x^l_{iq_i}$ from the corresponding support points.\\
\noindent  \hspace*{1cm} $b)$ If $\|c-s_j\|^2= \| c-s_l\|^2$ then\\ \noindent \hspace*{1.8cm} Identify the minimal mass $d_{\text{min}}= \min\limits_{i \leq N} d^l_{iq_i}$ among the $x^l_{iq_i}$.\\ \noindent \hspace*{1.8cm} Set  $d_l=d_l-d_{\text{min}}$ and $d^l_{iq_i}=d^l_{iq_i}-d_{\text{min}}$ for all $i \leq N$. \\ \noindent \hspace*{1.8cm} For all $i\leq N$, if $d^l_{iq_i}=0$, remove $x^l_{iq_i}$ from  $\supp(P_i^l)$ and reindex $P_i^l$ and $d^l_{i1},\dots,d^l_{i|P_i^l|}$.
\\ \noindent \hspace*{1.8cm} 
For all $i\leq N$, add $x^l_{iq_i}$ to $\supp(P_i^j)$ if it is not in it yet. In this case, $|P_i^j|$ increases\\ \noindent \hspace*{1.8cm}  by one and we index the support point as  $x^j_{i|P_i^j|}$ (with $d^j_{i|P_i^j|}=0$).  \\ \noindent \hspace*{1.8cm} Let now $p_i$ be such that $x^j_{ip_i}=x^l_{iq_i}$ for all $i \leq N$.\\ \noindent \hspace*{1.8cm} Set $d_j=d_j+d_{\text{min}}$ and $d^j_{ip_i}=d^j_{ip_i}+d_{\text{min}}$ for all $i \leq N$.\\ \noindent \hspace*{1.8cm} 
If $d_l>0$, go back to $a)$.
\item {\bf \normalsize  (Spread out each support point to a set of weighted centroids)}\\\small
\noindent For $l=1$ ascending to $l=r$ \\
\noindent  \hspace*{0.5cm} Create an empty partial measure $\bar{P^l}$.\\
\noindent  \hspace*{0.5cm} $a)$ For each $i\leq N$, identify the index $q_i$ for a lexicographically maximal support point  \\\noindent  \hspace*{0.9cm}  $x^l_{iq_i}$ in $P_i^l$. Then compute the weighted centroid $c=\sum_{i=1}^N \lambda_i x^l_{iq_i}$.\\
\noindent  \hspace*{0.5cm} $b)$ Identify the minimal mass $d_{\text{min}}= \min\limits_{i \leq N} d^l_{iq_i}$ among the $x^l_{iq_i}$.\\ \noindent \hspace*{0.9cm}
Set  $d_l=d_l-d_{\text{min}}$ and $d^l_{iq_i}=d^l_{iq_i}-d_{\text{min}}$ for all $i \leq N$.\\ \noindent \hspace*{0.9cm} 
For all $i\leq N$, if $d^l_{iq_i}=0$, remove $x^l_{iq_i}$ from  $\supp(P_i^l)$ and reindex $P_i^l$ and $d^l_{i1},\dots,d^l_{i|P_i^l|}$.\\ \noindent \hspace*{0.9cm} 
Add $c$ to $\supp(\bar{P^l})$ with mass $d_{\text{min}}$. \\ \noindent \hspace*{0.9cm} If $d_l>0$, go back to $a)$.
 \item {\bf \normalsize  (Combine a new measure)}
\small \\\noindent Combine the partial measures $\bar P^l$ to a measure $\bar{P}' = \sum_{l=1}^r \bar{P^l}$. Return $\bar P'$.
\end{enumerate}

  \caption{Recovery of non-mass split}
   \label{algo:nonmasssplit}
\end{algorithm}

\noindent Note $c$ satisfies $\|c-s_l\|^2\leq \|c-s_j\|^2$ (because of optimality of $\bar P_\text{org}$) and it is a maximizer of $\|c-s_l\|^2-\|c-s_j\|^2 \leq 0$. If $\|c-s_l\|^2=\|c-s_j\|^2$, which is checked in $2b)$, then mass can be shifted from $s_l$ to $s_j$ to make $(d_1,\dots,d_r)$ lexicographically larger, while keeping optimality. The remainder of $2b)$ is a technical description of this shift of mass. 

\noindent{\bf Step 3.} Next, we perform a (greedy) routine to spread out the mass of each $s_l$ to several support points. We do so by picking a set of lexicographically maximal support points $x^l_{iq_i}$ in each $P_i^l$ (i.e., we pick an $x^l_{iq_i}$ with a largest first coordinate, and among those one with a largest second coordinate, and so on). Then we move mass $d_{\text{min}}$ to the weighted centroid $c=\sum_{i=1}^N \lambda_i x^l_{iq_i}$, where $d_{\text{min}}$ is the minimal mass among the $d^l_{q_i}$. We repeat this scheme until the whole mass of $s_l$ has been spread out. The result is a partial measure $\bar P^l$ that has a non-mass splitting transport by construction. Then we continue with the next support point.

\noindent{\bf Step 4.} Finally, we combine the partial measures $\bar P^l$ from Step $3$ to a new measure. It is at least as good an approximation of an exact barycenter as $\bar P_\text{org}$: Step $1$ and $2$ do not change the cost of transport. In Step $3$, for any chosen set of support points $x^l_{iq_i}$ we put the corresponding mass on their weighted centroid, which is best-possible (at least as good as transport from $s_l$).

We sum up the favorable properties of the algorithm in Theorem \ref{thm:fixingtheproperties}. In addition to the existence of a non-mass splitting transport, and keeping a $2$-approximation error, we are able to bound the size of the support by the square of the bound in Proposition \ref{sparsethm}. We do {\em not} prove that the returned measure is an approximate barycenter (which implies optimality over the given support by definition). The associated non-mass splitting transport is trivial to construct, but we do {\em not} prove that this transport is optimal. Due to this, we have to be careful in the wording of the following statements (Theorems \ref{thm:fixingtheproperties} and \ref{cor:nonmasssplit}). A detailed proof is given in Appendix \ref{sec:localimprovement}.

\begin{theorem}\label{thm:fixingtheproperties}
Algorithm \ref{algo:nonmasssplit} returns a measure $\bar P'$ supported on a subset of $S$ with $\phi(\bar P') \leq 2 \cdot \phi(\bar P)$ and there is a non-mass splitting transport realizing this bound. Further  $|\bar P'|\leq (\sum_{i=1}^N |P_i|  - N + 1)^2$.
\end{theorem}

\begin{figure}[t]
\begin{center}
\subfloat[Measure $P_1$]{
\begin{tikzpicture}[scale=1.5]
\node (11) at (-1,0) {$ $};
\node (22) at (3,0) {$ $};

\node (11) at (0,1.5) {$\frac{1}{4}$};
\node (12) at (1,-0.5) {$\frac{1}{2}$};
\node (21) at (2,1.5) {$\frac{1}{4}$};

\fill [blue] (0,1) circle (3pt);
\fill [blue] (1,0) circle (4.5pt);
\fill [blue] (2,1) circle (3pt);
\fill [black] (0,0) circle (1pt);
\fill [black] (1,1) circle (1pt);
\fill [black] (2,0) circle (1pt);
\end{tikzpicture}
}
\qquad
\subfloat[Measure $P_2$]{
\begin{tikzpicture}[scale=1.5]
\node (11) at (-1,0) {$ $};
\node (22) at (3,0) {$ $};

\node (11) at (0,-0.5) {$\frac{1}{4}$};
\node (12) at (1,1.5) {$\frac{1}{2}$};
\node (21) at (2,-0.5) {$\frac{1}{4}$};

\fill [black] (0,1) circle (1pt);
\fill [black] (1,0) circle (1pt);
\fill [black] (2,1) circle (1pt);
\fill [red] (0,0) circle (3pt);
\fill [red] (1,1) circle (4.5pt);
\fill [red] (2,0) circle (3pt);
\end{tikzpicture}
}

\subfloat[Transport from $s_1,s_2$]{
\begin{tikzpicture}[scale=1.5]
\node (11) at (-1,0) {$ $};
\node (22) at (3,0) {$ $};

\fill [black] (0,1) circle (1pt);
\fill [black] (1,0) circle (1pt);
\fill [black] (2,1) circle (1pt);
\fill [black] (0,0) circle (1pt);
\fill [black] (1,1) circle (1pt);
\fill [black] (2,0) circle (1pt);

\fill [black] (0,0.5) circle (3pt);
\fill [black] (2,0.5) circle (5.25pt);

\fill [blue] (1,0) circle (4.5pt);
\fill [red] (1,1) circle (4.5pt);

\fill [blue] (0,1) circle (3pt);

\fill [blue] (2,1) circle (3pt);
\fill [red] (0,0) circle (3pt);

\fill [red] (2,0) circle (3pt);

\node (11) at (-0.25,.5) {$s_1$};
\node (21) at (2.35,0.5) {$s_2$};

\draw[black,line width=0.8mm,->] (0,0.5) -- (0,0.05);
\draw[black,line width=0.8mm,->] (2,0.5) -- (1.05,0.95);
\draw[black,line width=0.8mm,->] (2,0.5) -- (2,0.05);

\draw[black,line width=0.8mm,->] (0,0.5) -- (0,0.95);
\draw[black,line width=0.8mm,->] (2,0.5) -- (1.05,0.05);
\draw[black,line width=0.8mm,->] (2,0.5) -- (2,0.95);

\end{tikzpicture}}
\qquad 
\subfloat[A centroid $c$ (Step $2a)$)]{
\begin{tikzpicture}[scale=1.5]
\node (11) at (-1,0) {$ $};
\node (22) at (3,0) {$ $};

\fill [blue] (1,0) circle (4.5pt);
\fill [black] (1,0.5) circle (4.5pt);
\fill [red] (1,1) circle (4.5pt);

\fill [black] (0,1) circle (1pt);
\fill [black] (2,1) circle (1pt);
\fill [black] (0,0) circle (1pt);
\fill [black] (2,0) circle (1pt);

\fill [black] (0,0.5) circle (3pt);
\fill [black] (2,0.5) circle (3pt);

\node (21) at (1.25,0.65) {$c$};
\node (11) at (-0.25,.5) {$s_1$};
\node (21) at (2.25,0.5) {$s_2$};


\end{tikzpicture}} 

\subfloat[Mass shift to $s_1$ (Step $2b)$]{
\begin{tikzpicture}[scale=1.5]
\node (11) at (-1,0) {$ $};
\node (22) at (3,0) {$ $};
\fill [black] (0,1) circle (1pt);
\fill [black] (1,0) circle (1pt);
\fill [black] (2,1) circle (1pt);
\fill [black] (0,0) circle (1pt);
\fill [black] (1,1) circle (1pt);
\fill [black] (2,0) circle (1pt);

\fill [black] (0,0.5) circle (5.25pt);
\fill [black] (2,0.5) circle (3pt);

\fill [blue] (1,0) circle (4.5pt);
\fill [red] (1,1) circle (4.5pt);

\fill [blue] (0,1) circle (3pt);

\fill [blue] (2,1) circle (3pt);
\fill [red] (0,0) circle (3pt);

\fill [red] (2,0) circle (3pt);

\draw[black,line width=0.8mm,->] (0,0.5) -- (0,0.05);
\draw[black,line width=0.8mm,->] (0,0.5) -- (0.95,0.95);
\draw[black,line width=0.8mm,->] (2,0.5) -- (2,0.05);

\draw[black,line width=0.8mm,->] (0,0.5) -- (0,0.95);
\draw[black,line width=0.8mm,->] (0,0.5) -- (0.95,0.05);
\draw[black,line width=0.8mm,->] (2,0.5) -- (2,0.95);

\node (11) at (-0.35,.5) {$s_1$};
\node (21) at (2.25,0.5) {$s_2$};
\end{tikzpicture}}
\qquad 
\subfloat[Spread of $s_1$ (Step $3$)]{
\begin{tikzpicture}[scale=1.5]
\node (11) at (-1,0) {$ $};
\node (22) at (3,0) {$ $};

\fill [black] (0,1) circle (1pt);
\fill [black] (1,0) circle (1pt);
\fill [black] (2,1) circle (1pt);
\fill [black] (0,0) circle (1pt);
\fill [black] (1,1) circle (1pt);
\fill [black] (2,0) circle (1pt);

\fill [black] (0,0.5) circle (3pt);
\fill [black] (1,0.5) circle (4.5pt);

\fill [blue] (1,0) circle (4.5pt);
\fill [red] (1,1) circle (4.5pt);

\fill [blue] (0,1) circle (3pt);

\fill [red] (0,0) circle (3pt);

\draw[black,line width=0.8mm,->] (0,0.5) -- (0,0.05);
\draw[black,line width=0.8mm,->] (1,0.5) -- (1,0.95);

\draw[black,line width=0.8mm,->] (0,0.5) -- (0,0.95);
\draw[black,line width=0.8mm,->] (1,0.5) -- (1,0.05);

\end{tikzpicture}}

\subfloat[No spread of $s_2$ (Step $3)$]{
\begin{tikzpicture}[scale=1.5]
\node (11) at (-1,0) {$ $};
\node (22) at (3,0) {$ $};
\fill [black] (0,1) circle (1pt);
\fill [black] (1,0) circle (1pt);
\fill [black] (2,1) circle (1pt);
\fill [black] (0,0) circle (1pt);
\fill [black] (1,1) circle (1pt);
\fill [black] (2,0) circle (1pt);

\draw[black,line width=0.8mm,->] (2,0.5) -- (2,0.05);

\draw[black,line width=0.8mm,->] (2,0.5) -- (2,0.95);

\fill [black] (2,0.5) circle (3pt);



\fill [blue] (2,1) circle (3pt);

\fill [red] (2,0) circle (3pt);

\end{tikzpicture}}
\qquad 
\subfloat[$\bar P'$ (Step $4$)]{
\begin{tikzpicture}[scale=1.5]
\node (11) at (-1,0) {$ $};
\node (22) at (3,0) {$ $};

\node (11) at (0.25,.75) {$\frac{1}{4}$};
\node (12) at (1.25,0.75) {$\frac{1}{2}$};
\node (21) at (2.25,0.75) {$\frac{1}{4}$};

\fill [black] (0,0.5) circle (3pt);
\fill [black] (1,0.5) circle (4.5pt);
\fill [black] (2,0.5) circle (3pt);

\fill [black] (0,1) circle (1pt);
\fill [black] (1,0) circle (1pt);
\fill [black] (2,1) circle (1pt);
\fill [black] (0,0) circle (1pt);
\fill [black] (1,1) circle (1pt);
\fill [black] (2,0) circle (1pt);
\end{tikzpicture}}

  \end{center}
\caption{Two measures $P_1$, $P_2$ in the top row and a run of Steps $2 - 4$ of Algorithm \ref{algo:nonmasssplit} for given support points $s_1,s_2$ of mass $d_1=\frac{1}{4}$, $d_2=\frac{3}{4}$. Note $s_1,s_2 \notin S_{\text{org}}$, which may happen in later iterations of Algorithm \ref{algo:heuristic}, where Algorithm \ref{algo:nonmasssplit} is used as a subroutine.} 
\label{fig:exampleforstep2}
\end{figure}

Let us discuss a small example for Steps $2-4$ of the algorithm.

\begin{example}\label{ex0}
We revisit the measures $P_1$ and $P_2$ used for Example \ref{ex1} and again let $\lambda_1=\lambda_2=\frac{1}{2}$. 
They receive their mass transported from two fixed support points $s_1,s_2$ of mass $d_1=\frac{1}{4}$, $d_2=\frac{3}{4}$ (second row, left). Note that $s_1,s_2 \notin S_{\text{org}}$, which may happen in later iterations of Algorithm \ref{algo:heuristic} (Section \ref{sec:itersec}), where Algorithm \ref{algo:nonmasssplit} is used as a subroutine. (For this example, this does not matter.)

The two central points, which receive their mass from $s_2$, have a centroid $c$ that is equally far from $s_1$ and $s_2$ (second row, right). These two points would be selected in Step $2a)$ of Algorithm \ref{algo:nonmasssplit} and their mass shifted from $s_2$ to $s_1$ in Step $2b)$. Then $d_1=\frac{3}{4}$, $d_2=\frac{1}{4}$ (third row, left).

In Step $3$, the mass of $s_1$ and $s_2$ is spread out to a set of centroids that transport to just a single support point in each measure. The result for $s_1$ is depicted in the third row (right). By lexicographically maximal choice of the points in the measures, the central support point of mass $\frac{1}{2}$ is constructed first, followed by the left one of mass $\frac{1}{4}$. $s_2$ is not changed, because it already is the centroid of a set of single support points in each measure (fourth row, left).

These partial measures are combined to form measure $\bar P'$ in Step $4$ (fourth row, right) and the algorithm stops. We actually found an exact barycenter, which is not the case in general.\hfill$\square$ 
\end{example}

We close our discussion of Algorithm \ref{algo:nonmasssplit} by showing that it runs in strongly-polynomial time. The quite technical proof is given in Appendix \ref{sec:localimprovement}.

\begin{theorem}\label{cor:nonmasssplit}
For all rational input, a measure can be computed in strongly-polynomial time that is a $2$-approximation of a barycenter and for which there is a non-mass splitting transport realizing this bound.
\end{theorem}

\subsection{An Iterative Improvement}\label{sec:itersec}

Finally, we combine Algorithms \ref{algo:2approx} and \ref{algo:nonmasssplit} to an iterative scheme, denoted as Algorithm \ref{algo:heuristic}. The algorithm begins by computing an approximate barycenter in $S_{\text{org}}$ using Algorithm \ref{algo:2approx}. Then Algorithm \ref{algo:nonmasssplit} is used to spread out its support points to find a new measure $\bar P'$ of better approximation error and that allows for a non-mass splitting transport. We set $S_{\text{org}}=\text{supp}(\bar P')$ and repeat Algorithm \ref{algo:2approx} to find an optimal approximate barycenter over this new support (in other words, an optimal vertex of LP (\ref{baryLP}) over the new support is found). Then its support points are spread out again using Algorithm \ref{algo:nonmasssplit}. This scheme is repeated until there is no improvement anymore. 

After a finite number of iterations, the algorithm terminates with a sparse $2$-approximate barycenter supported on a subset of $S$, and with a non-mass splitting optimal transport. This is a provable approximation that possesses both favorable properties of an exact barycenter, sparsity and non-mass split, at the same time.

\begin{algorithm}[H]
\vspace{0.2cm}

 {\bf Input}
\begin{itemize}
\item Measures $P_1,\dots,P_N \subset \mathbb{R}^d$
\item $\lambda_1,\dots,\lambda_N> 0$ with $\sum_{i=1}^N \lambda_i =1$
\end{itemize}
 {\bf Algorithm}
\begin{enumerate}
\item Compute a (sparse) $2$-approximate barycenter $\bar P_\text{org}$ in $S_\text{org}$ (and an optimal transport) \\
\noindent using Algorithm \ref{algo:2approx}.
\item Use $\bar P_\text{org}$ (and its transport) as input for Algorithm \ref{algo:nonmasssplit} to find a measure $\bar P'$.\\ \noindent    
If $\bar P'\neq \bar P_\text{org}$, set $S_\text{org}=\text{supp}(\bar P')$ and go back to $1$. Else return $\bar P'$.
\end{enumerate}
  \caption{Iterative improvement}
   \label{algo:heuristic}
\end{algorithm}

\begin{theorem}\label{thm:heuristic}
Algorithm \ref{algo:heuristic} returns an approximate barycenter $\bar P'$ supported on a subset of $S$ for which $\phi(\bar P') \leq 2 \cdot \phi(\bar P)$, where  $\bar P$ is a barycenter, and there is a non-mass splitting optimal transport realizing this bound. Further  $|\bar P'|\leq \sum_{i=1}^N |P_i|  - N + 1$.
\end{theorem}

We prove Theorem \ref{thm:heuristic} in Appendix \ref{sec:iterativeproofs}. Let us take a closer look at the approximation error of Algorithm \ref{algo:heuristic}. We distinguish three different measures:  $\bar P$ is an exact barycenter, $\bar P_{\text{org}}$ is an approximate barycenter in $\mathcal{P}_{\text{org}}^2(\R^d)$ , and $\bar P'$ is the solution of Algorithm \ref{algo:heuristic}. By optimality of $\bar P$ and $\bar P_{\text{org}}$ with respect to $\phi$ in their respective support, we have
$$\phi(\bar P)\leq \phi(\bar P') \leq\phi(\bar P_{\text{org}}).$$
We are particularly interested in the gap between $\phi(\bar P)$ and $\phi(\bar P')$. Theorem \ref{thm:originalsupport} states $\phi(\bar P_{\text{org}}) \leq 2 \cdot \phi(\bar P) $. Thus the whole sequence of inequalities is bounded by a total approximation factor of $2$. This implies that if $\alpha \phi(\bar P')=\phi(\bar P_{\text{org}})$ for some $\alpha \geq 1$, then $\phi(\bar P') \leq \frac{2}{\alpha} \phi(\bar P)$. Informally, Algorithm \ref{algo:heuristic} already begins with a provable $2$-approximation and any improvement towards $\bar P'$ allows for the statement of a better approximation guarantee.

In practice, one obtains a strictly better approximation factor than $2$ for essentially all real-world problems using Algorithm \ref{algo:heuristic}. But there exist worst-case examples, such as the following, that show the bound is tight.

\begin{example} Consider the example depicted in Figure \ref{fig:factor2error}. Four measures $P_1,...,P_4$ are shown in the top row, $P_2$ and $P_3$ are depicted in the center. Note $P_2=P_3$. Each of the measures consists of two support points of mass $\frac{1}{2}$. Let $\epsilon >1$. $P_1$ is supported on coordinates $(-\epsilon,0)$ and $(\epsilon,1)$, $P_2$ and $P_3$ are supported on $(0,0)$ and $(0,1)$, and $P_4$ is supported on $(-\epsilon,1)$ and $(\epsilon,0)$, where $\epsilon > 1$. For increasing $\epsilon$, the horizontal distance of the support points of $P_1$ and $P_4$ to those of $P_2,P_3$ increases proportionally (second row).

Let $\lambda_{i}=\frac{1}{4}$ for $i=1,\dots,4$. Independently of $\epsilon$, an approximate barycenter $\bar P_{\text{org}}$ in $S_{\text{org}}$ is identical to $P_2=P_3$,  (third row, left). A corresponding optimal transport sends the mass to the support points in the same `layer' (third row, middle). Note that the support points are already the (weighted) centroids of the points they transport to, and that the transport is non-mass splitting. Because of this, Algorithm \ref{algo:heuristic} stops without any change to $\bar P_{\text{org}}$ at the end of the first iteration. 

The cost of transport for $\bar P_{\text{org}}$ is $\phi(\bar P_{\text{org}})=\frac{1}{4} \cdot 2\epsilon^2= \frac{1}{2}\epsilon^2$. (Recall $\lambda_{i}=\frac{1}{4}$ for all $i$.) An exact barycenter $\bar P$ (third row, right) and a corresponding optimal transport (fourth row) are strictly better. The coordinates for the two support points are $(-\frac{1}{2}\epsilon, \frac{3}{4})$ and $(\frac{1}{2}\epsilon, \frac{1}{4})$. The cost of transport is $\phi(\bar P)=\frac{1}{4}\cdot (\frac{3}{4}+ \epsilon^2)=\frac{3}{16} + \frac{1}{4}\epsilon^2$. For $\epsilon \rightarrow \infty$, 
$$ \frac{\phi(\bar P_{\text{org}})}{\phi(\bar P)}=\frac{\frac{1}{2}\epsilon^2}{\frac{3}{16}+ \frac{1}{4}\epsilon^2}\rightarrow 2.$$
Thus the error bound goes to $2$. \hfill$\square$ 
\end{example}

\begin{figure}
\begin{center}
\subfloat[Measure $P_1$]{
\begin{tikzpicture}[scale=1.5]
\node (11) at (0,-0.5) {$\frac{1}{2}$};
\node (21) at (2,1.5) {$\frac{1}{2}$};

\fill [blue] (0,0) circle (4.5pt);
\fill [black] (0,1) circle (1pt);
\fill [black] (1,0) circle (1pt);
\fill [black] (1,1) circle (1pt);
\fill [blue] (2,1) circle (4.5pt);
\fill [black] (2,0) circle (1pt);
\end{tikzpicture}}
\qquad \qquad  \qquad 
\subfloat[Measures $P_2,P_3$]{
\begin{tikzpicture}[scale=1.5]
\node (11) at (1,-0.5) {$\frac{1}{2}$};
\node (12) at (1,1.5) {$\frac{1}{2}$};

\fill [green] (1,0) circle (4.5pt);
\fill [green] (1,1) circle (4.5pt);

\fill [black] (0,1) circle (1pt);
\fill [black] (2,1) circle (1pt);
\fill [black] (0,0) circle (1pt);
\fill [black] (2,0) circle (1pt);
\end{tikzpicture}}
\qquad \qquad  \qquad 
\subfloat[Measure $P_4$]{
\begin{tikzpicture}[scale=1.5]
\node (11) at (0,1.5) {$\frac{1}{2}$};
\node (21) at (2,-0.5) {$\frac{1}{2}$};

\fill [black] (0,0) circle (1pt);
\fill [red] (0,1) circle (4.5pt);
\fill [black] (1,0) circle (1pt);
\fill [black] (1,1) circle (1pt);
\fill [black] (2,1) circle (1pt);
\fill [red] (2,0) circle (4.5pt);
\end{tikzpicture}
}

\subfloat[all Measures, $\epsilon\to\infty$]{
\begin{tikzpicture}[scale=1.5]

\node (33) at (1,-0.4) {$ $};
\node (44) at (1,1.3) {$ $};

\draw[black,line width=0.6mm,dotted,->] (2,0) -- (5,0);
\draw[black,line width=0.6mm,dotted,->] (0,0) -- (-3,0);

\draw[black,line width=0.6mm,dotted,->] (2,1) -- (5,1);
\draw[black,line width=0.6mm,dotted,->] (0,1) -- (-3,1);

\fill [green] (1,0) circle (4.5pt);
\fill [green] (1,1) circle (4.5pt);

\fill [red] (-1,1) circle (4.5pt);
\fill [red] (3,0) circle (4.5pt);

\fill [blue] (-1,0) circle (4.5pt);
\fill [blue] (3,1) circle (4.5pt);

\end{tikzpicture}}

\subfloat[Measure $\bar P_{\text{org}}$]{
\begin{tikzpicture}[scale=1.5]
\node (11) at (1,-0.5) {$\frac{1}{2}$};
\node (12) at (1,1.5) {$\frac{1}{2}$};

\fill [black] (1,0) circle (4.5pt);
\fill [black] (1,1) circle (4.5pt);

\fill [black] (0,1) circle (1pt);
\fill [black] (1,0) circle (1pt);
\fill [black] (2,1) circle (1pt);
\fill [black] (0,0) circle (1pt);
\fill [black] (1,1) circle (1pt);
\fill [black] (2,0) circle (1pt);
\end{tikzpicture}}
\qquad \qquad  \qquad
\subfloat[Transport for $\bar P_{\text{org}}$]{
\begin{tikzpicture}[scale=1.5]

\fill [blue] (0,0) circle (4.5pt);

\fill [blue] (2,1) circle (4.5pt);
\fill [red] (0,1) circle (4.5pt);
\fill [red] (2,0) circle (4.5pt);

\draw[black,line width=0.6mm,->] (1,0) -- (0.05,0);
\draw[black,line width=0.6mm,->] (1,0) -- (1.95,0);

\draw[black,line width=0.6mm,->] (1,1) -- (0.05,1);
\draw[black,line width=0.6mm,->] (1,1) -- (1.95,1);

\fill [green] (1,0) circle (4.5pt);
\fill [green] (1,1) circle (4.5pt);

 \path[black,line width=0.7mm,->] (1,1.1)   edge[out=135,in=45, my loop] node  {} (1,1.1);
 \path[black,line width=0.7mm,->] (1,-0.1)   edge[out=235,in=315, my loop] node  {} (1,-0.1);

\end{tikzpicture}}
\qquad \qquad  \qquad 
\subfloat[Barycenter $\bar P$]{
\begin{tikzpicture}[scale=1.5]
\node (33) at (1,-0.4) {$ $};

\node (44) at (1,1.3) {$ $};

\node (11) at (1.5,-0.25) {$\frac{1}{2}$};
\node (12) at (0.5,1.25) {$\frac{1}{2}$};

\fill [black] (0.5,0.75) circle (4.5pt);
\fill [black] (1.5,0.25) circle (4.5pt);

\fill [black] (0,1) circle (1pt);
\fill [black] (1,0) circle (1pt);
\fill [black] (2,1) circle (1pt);
\fill [black] (0,0) circle (1pt);
\fill [black] (1,1) circle (1pt);
\fill [black] (2,0) circle (1pt);
\end{tikzpicture}}

\subfloat[Transport for $\bar P$, $\epsilon\to\infty$]{
\begin{tikzpicture}[scale=1.5]

\node (33) at (1,-0.4) {$ $};
\node (44) at (1,1.3) {$ $};

\draw[black,line width=0.6mm,dotted,->] (2,0) -- (5,0);
\draw[black,line width=0.6mm,dotted,->] (0,0) -- (-3,0);

\draw[black,line width=0.6mm,dotted,->] (2,1) -- (5,1);
\draw[black,line width=0.6mm,dotted,->] (0,1) -- (-3,1);

\fill [green] (1,0) circle (4.5pt);
\fill [green] (1,1) circle (4.5pt);

\fill [red] (-1,1) circle (4.5pt);
\fill [red] (3,0) circle (4.5pt);

\fill [blue] (-1,0) circle (4.5pt);
\fill [blue] (3,1) circle (4.5pt);

\fill [black] (0,0.75) circle (4.5pt);
\fill [black] (2,0.25) circle (4.5pt);

\draw[black,line width=0.6mm,->]  (0,0.75) -- (-0.95,0);
\draw[black,line width=0.6mm,->]  (0,0.75) -- (-0.95,0.95);
\draw[black,line width=0.6mm,->]  (0,0.75) -- (0.95,0.95);

\draw[black,line width=0.6mm,->] (2,0.25) -- (2.95,0.05);
\draw[black,line width=0.6mm,->] (2,0.25) -- (1.05,0.05);
\draw[black,line width=0.6mm,->] (2,0.25) -- (2.95,0.95);

\end{tikzpicture}}

  \end{center}
\caption{Four measures $P_1,...,P_4$ (depicted for $\epsilon=1$) in the first row. Note $P_2=P_3$. For increasing $\epsilon$, the horizontal distance of the support of $P_1$ and $P_4$ to $P_2=P_3$ increases (second row). An approximate barycenter $\bar P_{\text{org}}$ in $S_{\text{org}}$, corresponding transport, and an exact barycenter $\bar P$ (all depicted for $\epsilon=1$) in the third row. The transport for $\bar P$ in the fourth row. Algorithm \ref{algo:heuristic} returns $\bar P_{\text{org}}$. For $\epsilon\to\infty$, $\frac{\phi(\bar P_{\text{org}})}{\phi(\bar P)}\rightarrow 2$, i.e., the error goes to $2$.}
\label{fig:factor2error}
\end{figure}
 
In Section \ref{computations} we conclude the paper with a discussion of the theoretical scaling of our algorithms and some observations on practical computations. In our implementation, we use some tweaks for a speed-up of Algorithm \ref{algo:heuristic}. First, we perform Step $3$ of Algorithm \ref{algo:nonmasssplit} as the exact computation of a barycenter $\bar P^l$ when the number $N_l$ of support points to which a given $s_l$ transports is low. This leads to a better approximation bound for $\bar P'$ at the end of each iteration of Algorithm \ref{algo:heuristic} and a lower number of iterations overall. Further, this leads to a smaller support for the LPs in the second iteration and beyond.

Recall $N_l$ is bounded below by $N$, as the support point transports mass to at least one support point in each measure. We tried different values for $N_l$ for the cutoff to an exact barycenter computation $\bar P^l$. There is a tradeoff between each run of Algorithm \ref{algo:nonmasssplit} taking longer and a reduction in the total number of iterations. We observed good results for $N_l$ between $2\cdot N$ and $4\cdot N$ when $N\ll |S_{\text{org}}|$ (common for grid-structured data) and $N_l$ between $N+\log{N}$ and $N+4\log{N}$ when $N\gg |S_{\text{org}}|$ (common for data in general position). For these values, each run of Algorithm \ref{algo:nonmasssplit}, respectively Step $2$ of Algorithm \ref{algo:heuristic}, takes slightly longer, but the running time remains negligible compared to Step $1$ in each iteration. For the computations in Section \ref{computations}, we chose cutoffs of $2\cdot N$ and $N+\log{N}$. We observed a noticable positive impact, dropping the total running time of Algorithm  \ref{algo:heuristic} by about $20\%$ on average.


Second, we explicitly construct the associated transport devised in Algorithm \ref{algo:nonmasssplit} and use it for a warm-start of the subsequent Step $1$ of Algorithm \ref{algo:heuristic}. The LPs of Step $1$ are solved through a primal simplex method. By construction, the transport from the previous iteration is not only a feasible vertex of the new primal LP, but already close to the new optimum. (It is a set of weighted centroids transporting to their respective support points, after all.) Thus, the subsequent primal LP can be warm-started and finding the exact optimum over the new support is much faster than solving the LP from scratch. This is a crucial part of the implementation, as the LPs in later iterations can have millions of variables and otherwise would be slow to solve. With this tweak, the setup of the LP in iteration $2$ remains as a bottleneck. The LPs in later iterations can be set up through an update of the previous one. Because of the warm-start, the actual solution time of the LPs is negligible in comparison to the setup time.

\section{Sample Computations and Scaling}\label{computations}

We implemented Algorithms \ref{algo:2approx}, \ref{algo:nonmasssplit}, and \ref{algo:heuristic} in the Julia language using Clp as linear programming solver. Julia is a modern programming language for high-performance numerical computing that provides a competitive tradeoff between efficient, but cumbersome low-level languages (C, C\texttt{++}) and easy-to-use, but typically slow high-level languages (Python, Matlab) \citep{beks-14,LubinDunningIJOC}. A primal simplex method is called in Clp for the availability to warm-start iterations. The algorithms were run on a standard laptop (Win 10, 32GB memory, i7-6820). 

Our sample computations are on two representative types of data: the MNIST database of handwritten digits (widely-used for benchmarking) \citep{lbh-98} for grid-structured data, as well as the firehouse example from \cite{abm-16} (and randomly generated larger instances) for data in general position, i.e., where $S$ would be exponentially-sized. Together, these two settings cover most applications in practice: grid-structured data is common in machine learning and there is a wealth of algorithms for this setting. Data in general position typically arises when geographical locations are involved, like in many applications of operations research. The algorithms in this paper work on any data, and for any choice of $\lambda$. In contrast, working with data in general position, differing masses on the support points, or a non-uniform weight vector makes most of the algorithms in the literature impractical. Our goal is to identify for which data the practical performance of our algorithms is the most favorable. As we will see, it is data without an underlying grid-structure, with a small support and a large number of measures, that forms a best-case scenario. (We also treat the less favorable, but common grid setting in detail, for the sake of completeness.) We start with some sample runs for both types of data and then turn to the theoretical and practical scaling of computations.

\subsection{Sample runs of Algorithm \ref{algo:heuristic} and Observations}\label{sec:sampleruns}

\begin{figure}
\hspace*{-0.5cm}
\subfloat[Measure $P_1$]{
  \includegraphics[scale=0.28]{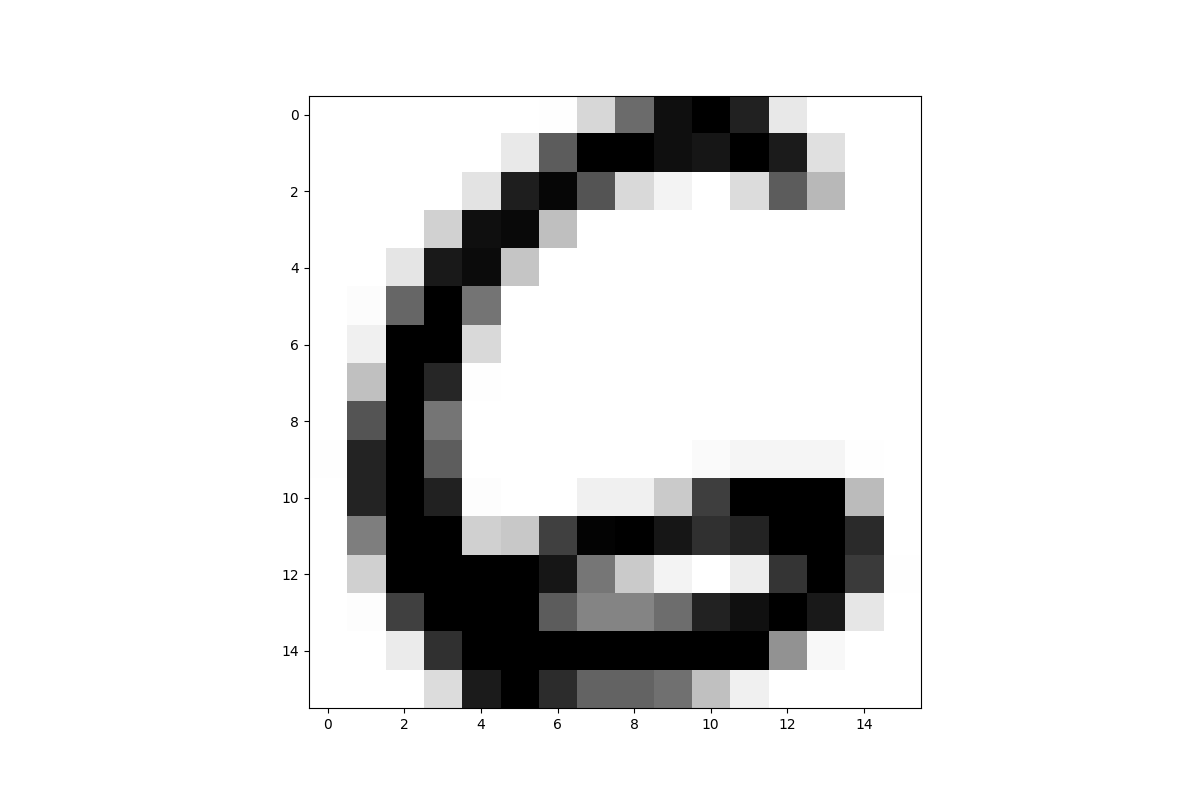}
}\hspace*{-2cm}
\subfloat[Measure $P_2$]{
  \includegraphics[scale=0.28]{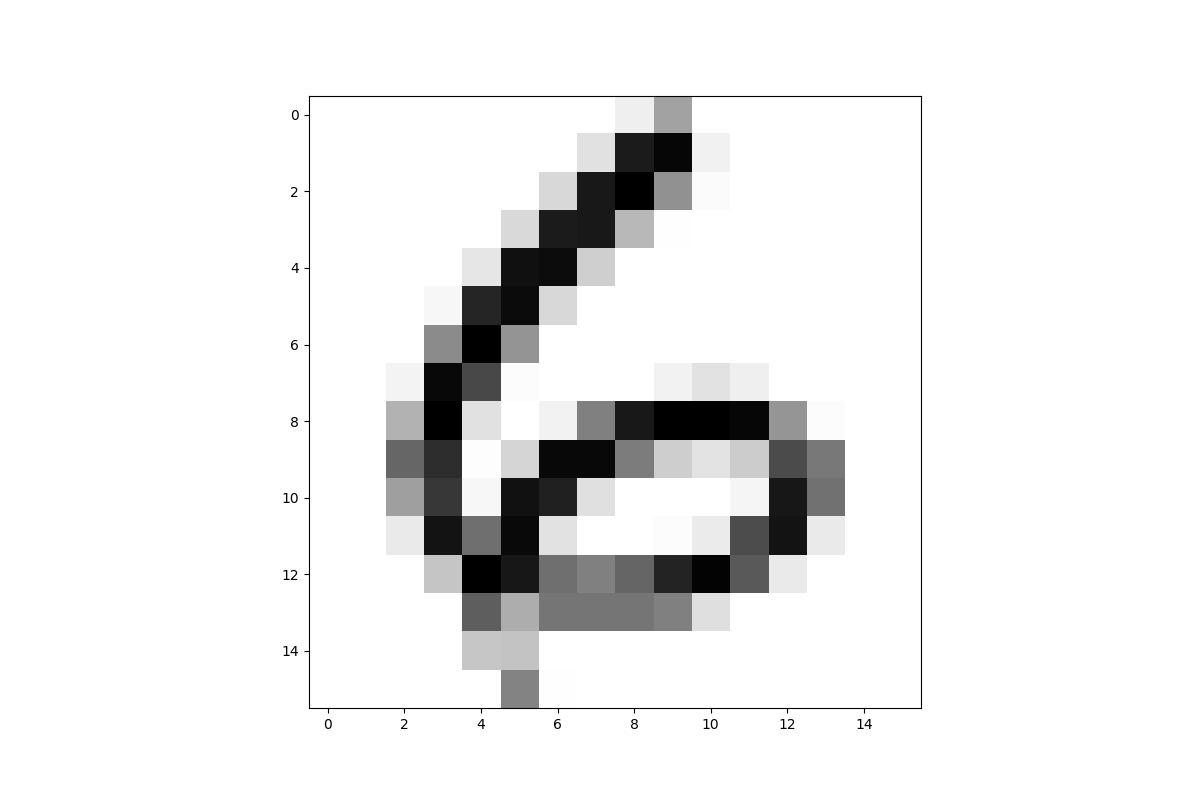}
}

\hspace*{-0.5cm}
\subfloat[Measure $P_3$]{
  \includegraphics[scale=0.28]{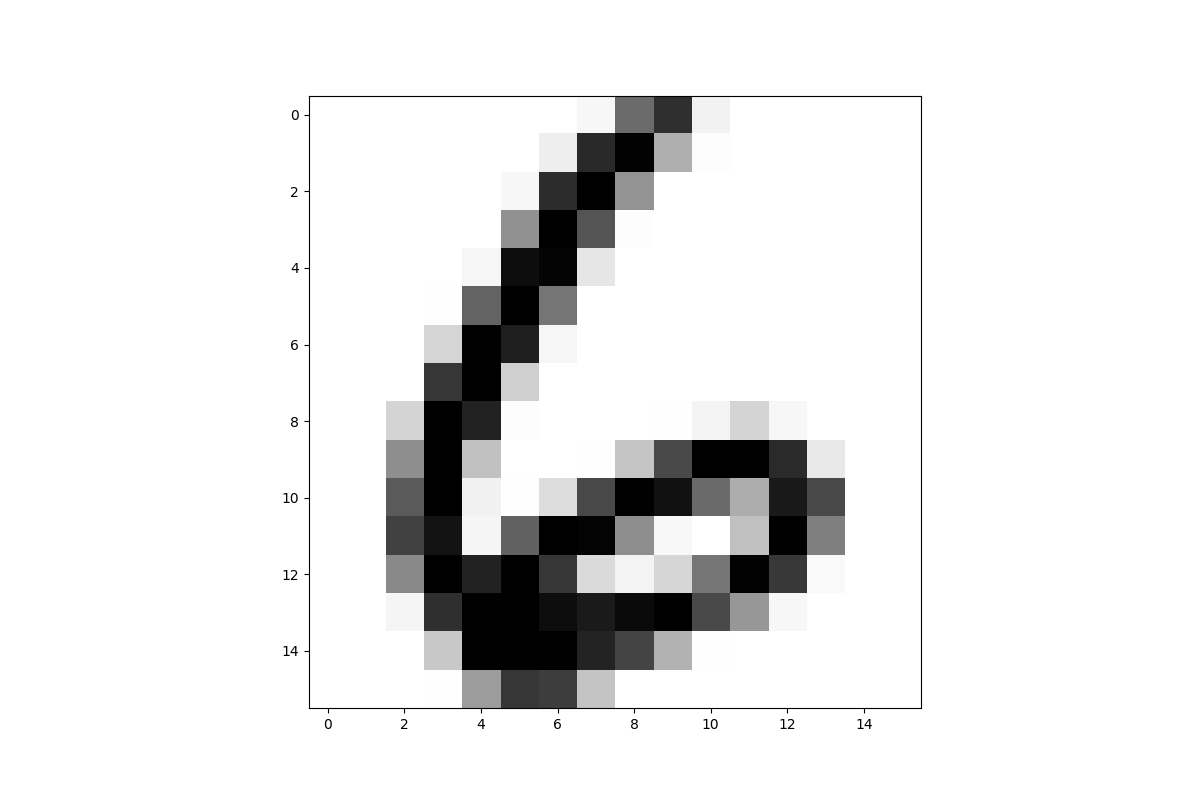}
}\hspace*{-2cm}
\subfloat[Measure $P_4$]{
  \includegraphics[scale=0.28]{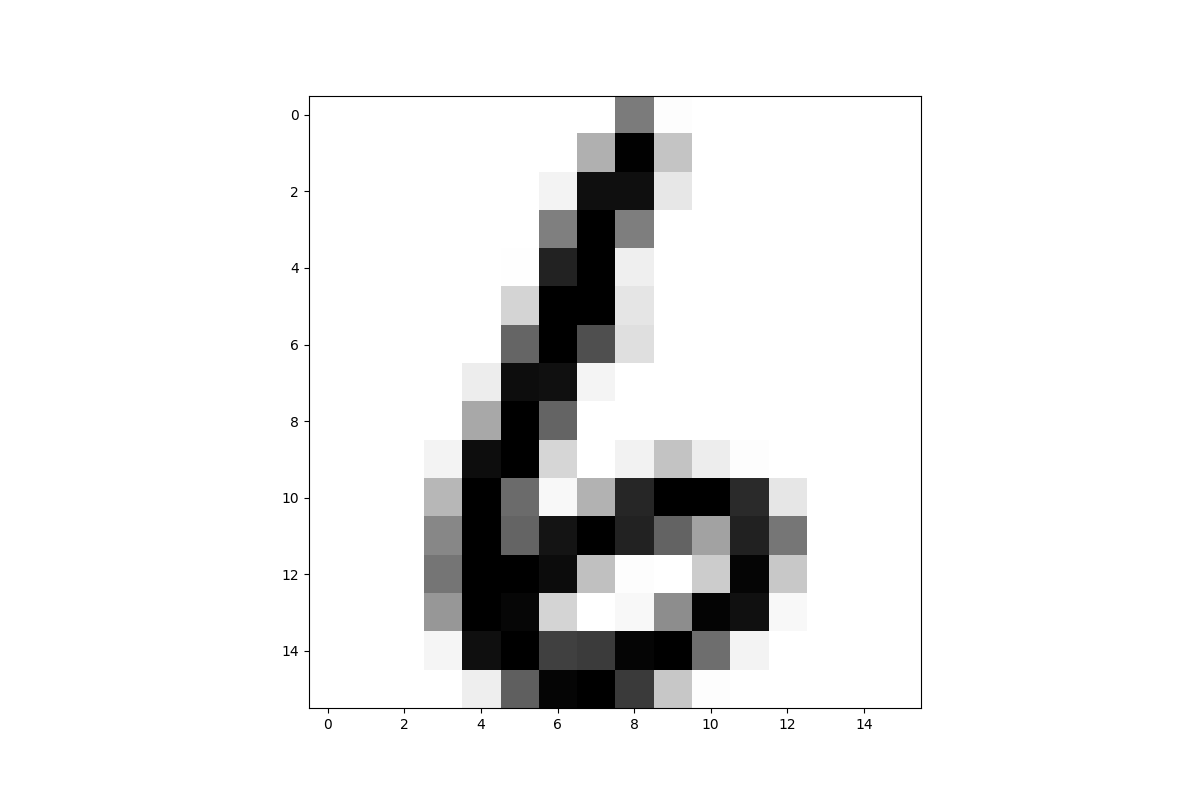}
}

\hspace*{+3cm}
\subfloat[Barycenter $\bar P$]{
  \includegraphics[scale=0.28]{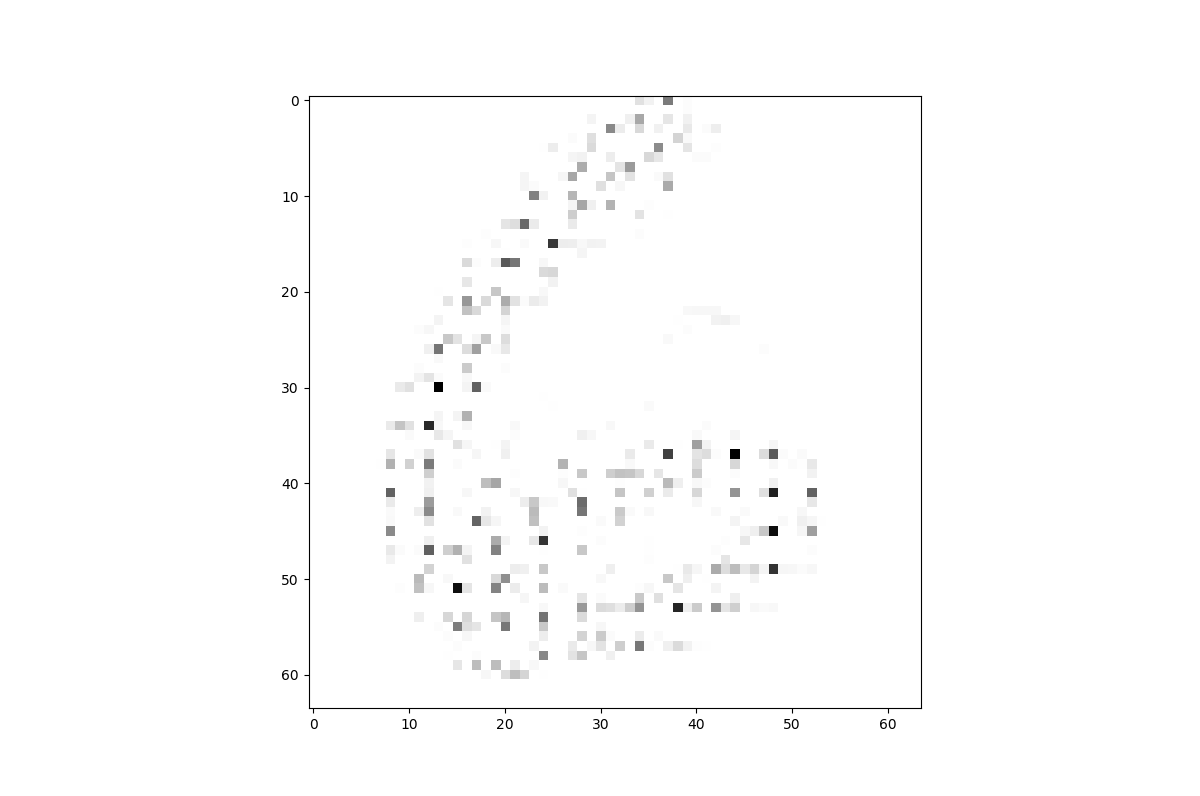}
}

\caption{Four measures $P_1,\dots,P_4$, scans of handwritten digits six, supported on a $16 \times 16$ grid. The barycenter $\bar P$ at the bottom.} 
\label{fig:foursixes}
\end{figure}

We begin by performing a sample run of Algorithm \ref{algo:heuristic} in a grid setting, using the four digits representing number six in a $16\times 16$ grid depicted in Figure \ref{fig:foursixes}. They have a barycenter depicted at the bottom of the figure (for all $\lambda_i=\frac{1}{4}$).

\begin{figure}[t]

\hspace*{-0.5cm}
\subfloat[Iteration $1$, Step $1$, Error $14.2\%$]{
  \includegraphics[scale=0.28]{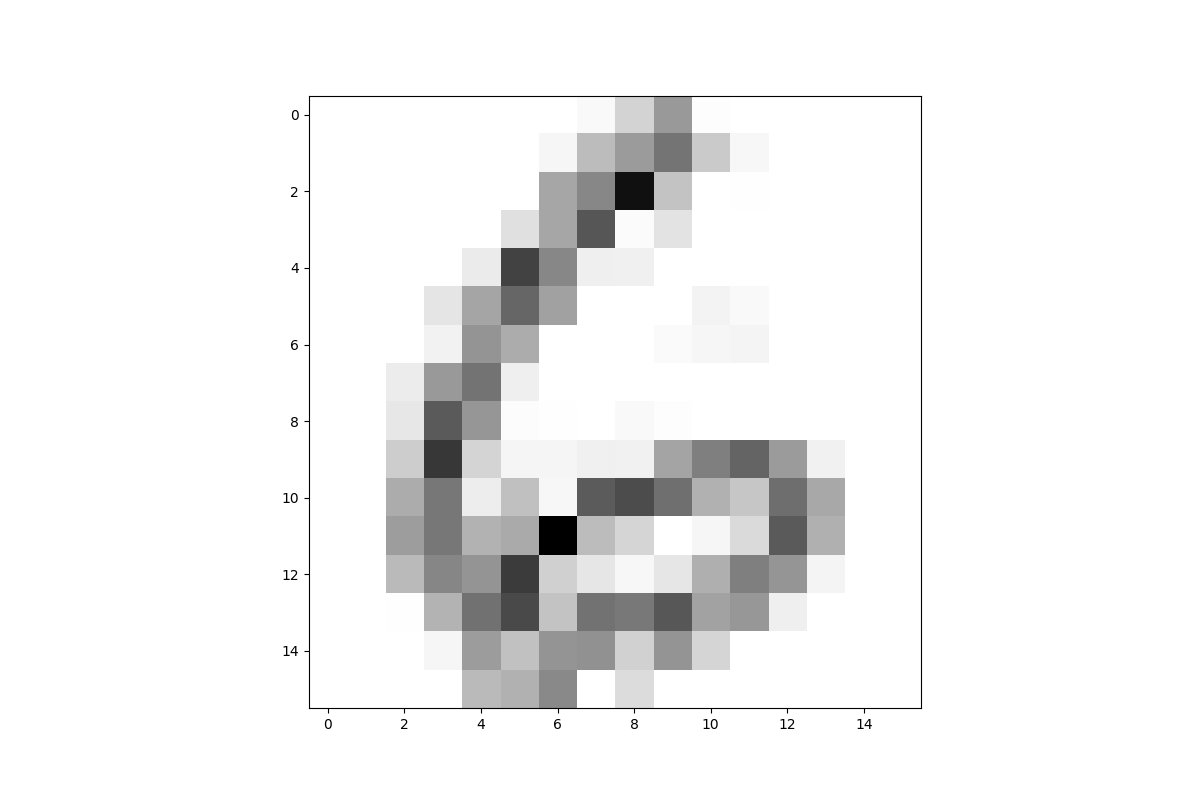}
}\hspace*{-2cm}
\subfloat[Iteration $1$, Step $2$, Error $4.3\%$]{
  \includegraphics[scale=0.28]{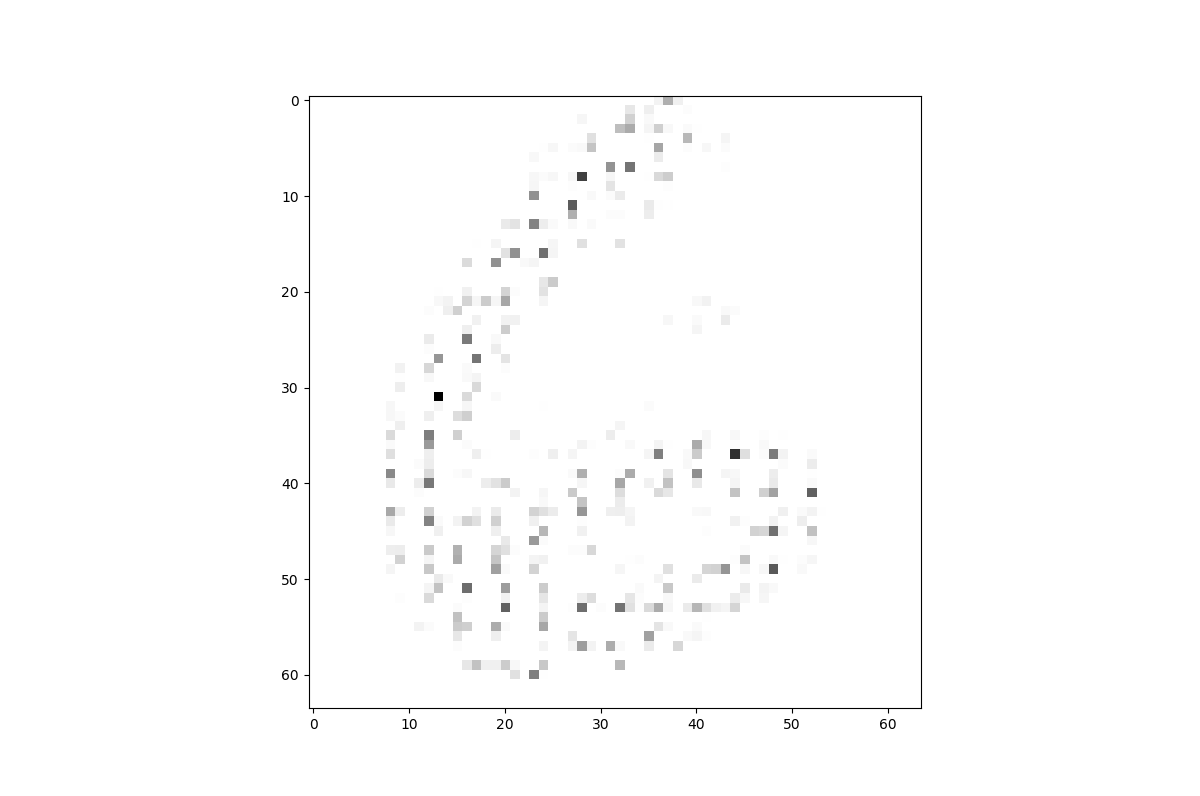}
}

\hspace*{-0.5cm}
\subfloat[Iteration $2$, Step $1$, Error $2.0\%$]{
  \includegraphics[scale=0.28]{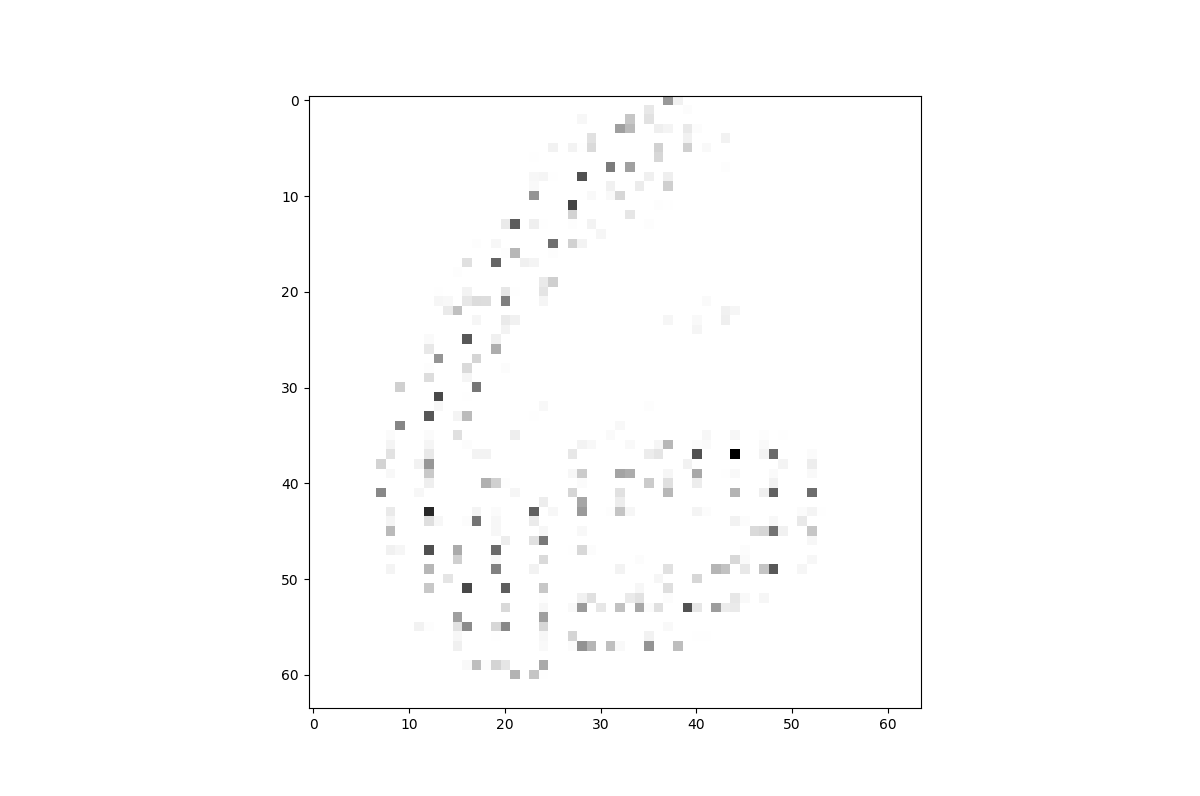}
}\hspace*{-2cm}
\subfloat[Iteration $2$, Step $2$, termination]{
  \includegraphics[scale=0.28]{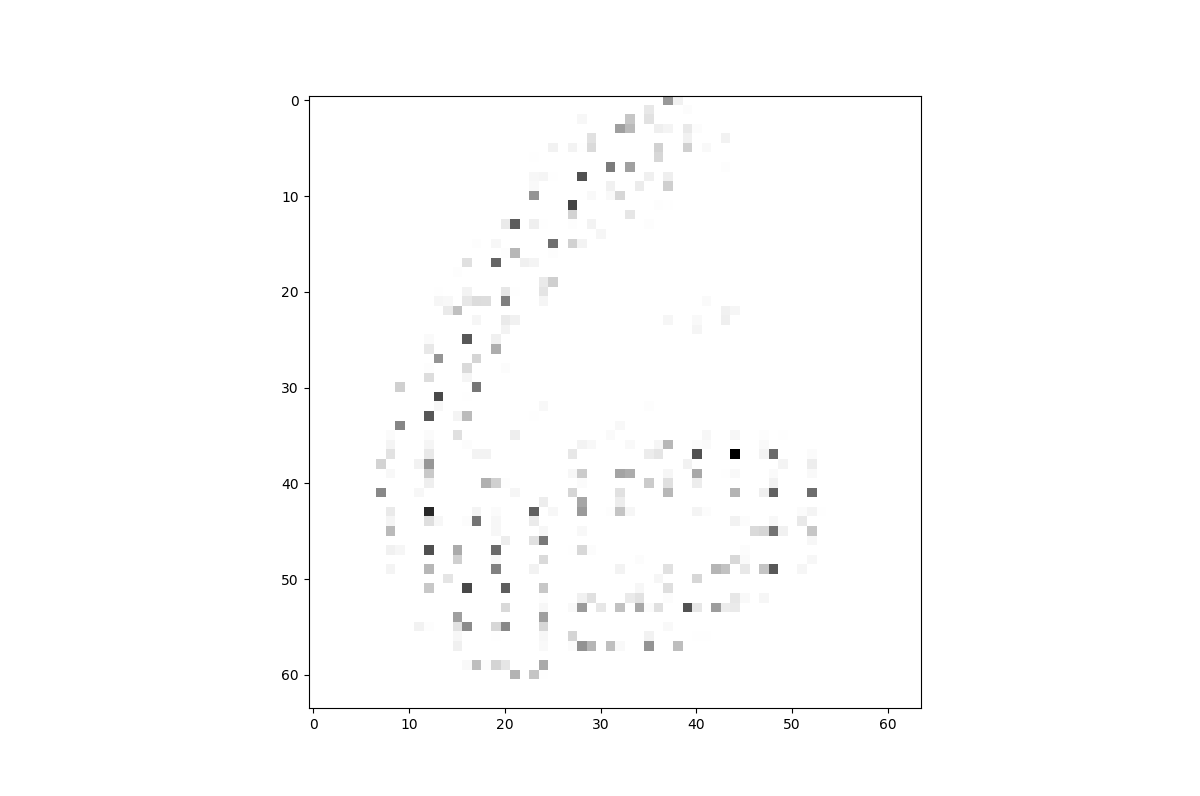}
}

\caption{A sample run of Algorithm \ref{algo:heuristic} for the measures in Figure \ref{fig:foursixes}. It already terminates after $2$ iterations.} 
\label{fig:stages}
\end{figure}

Figure \ref{fig:stages} shows the stages of a run of Algorithm \ref{algo:heuristic} for this input. Each row shows one of the iterations. The approximate barycenter in the original support is already a $1.142$--approximation of the exact barycenter (top left), i.e., $\phi(\bar P_\text{org})\leq 1.142 \cdot \phi(\bar P)$, which we denote as an additive $14.2\%$-error in the figure. The first split-up using Algorithm \ref{algo:nonmasssplit} (Steps $2$ to $4$) gives an improvement to a $4.3\%$-error (top right). This is further improved to a $2.0\%$-error (in Step $1$ of Iteration $2$) by computing an optimum over the support of the previous approximation (bottom left). Now the algorithm terminates, because all of the support points of this approximate barycenter are already the weighted centroids of the support points to which they transport mass, and there is no mass split. Algorithm \ref{algo:heuristic} completes in about $10$ seconds on average for a set of four measures ($9.6$ seconds for the above example). In contrast, the computation of an exact barycenter takes roughly $120$ seconds.

\begin{figure}
\begin{center}
\hspace*{-0.5cm}
\subfloat[Measure $P_1$]{
  \includegraphics[scale=0.76]{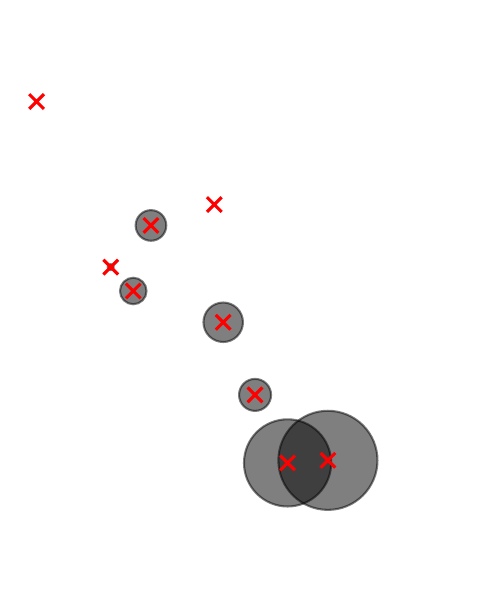}
}\hspace*{-1cm}
\subfloat[Measure $P_2$]{
  \includegraphics[scale=0.76]{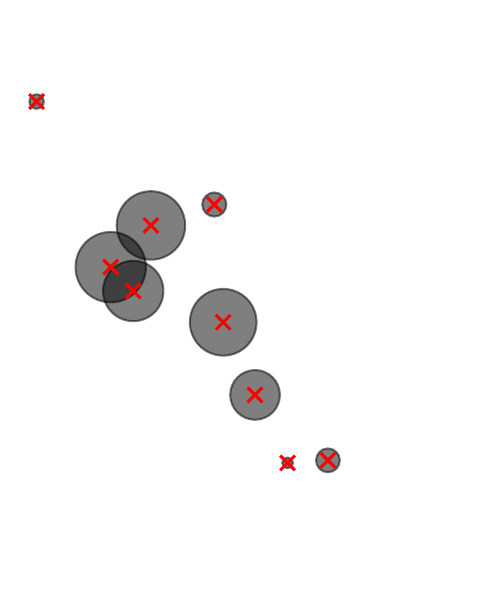}
}\hspace*{-1cm}
\subfloat[Barycenter $\bar P$]{
  \includegraphics[scale=0.76]{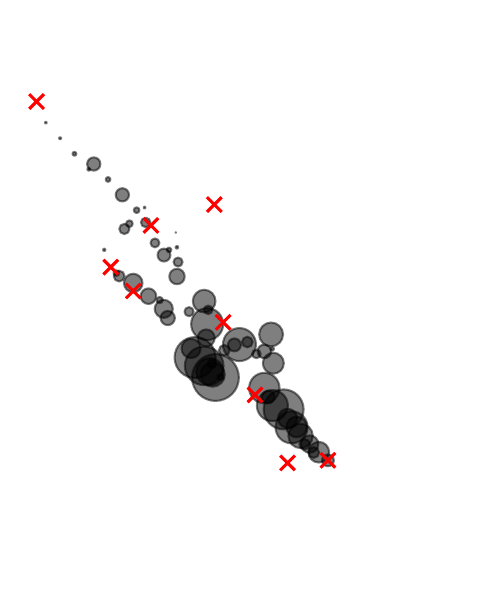}
}
\end{center}
\caption{Two (of eight) measures from a data set where the support points do not lie on a grid. All measures have the same support points with varying masses. The barycenter $\bar P$ to the right.} 
\label{fig:twogeneral}
\end{figure}

Next, we perform a sample run for data in general position. Here we begin with the input for the computations in \cite{abm-16}. There are $8$ measures with the same $9$ support points of varying masses. Figure \ref{fig:twogeneral} shows two of the measures. Circles of larger radius indicate higher mass. 

\begin{figure}
\begin{center}
\subfloat[$\bar P_{\text{org}}$, Error $10.2\%$]{
  \includegraphics[scale=0.76]{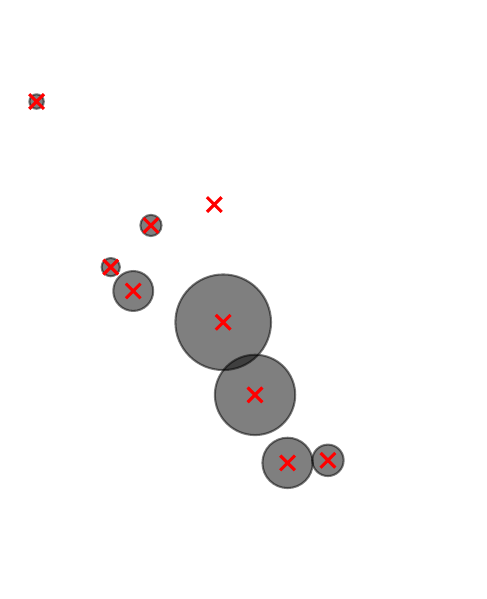}
}
\subfloat[first $\bar P'$, Error $1.9\%$]{
  \includegraphics[scale=0.76]{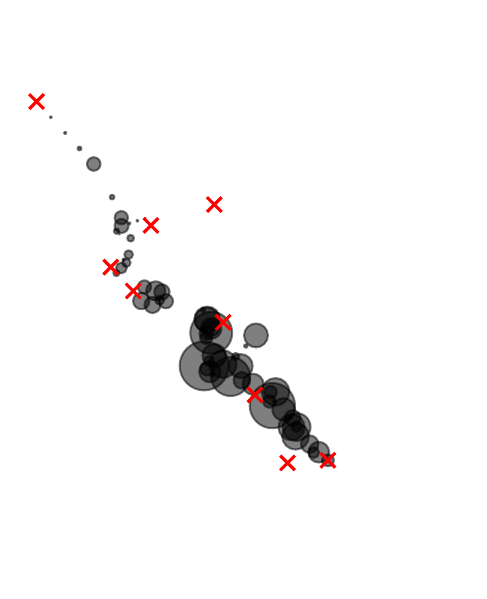}
}
\end{center}
\caption{Measures $\bar P_{\text{org}}$ and $\bar P'$ from the first iteration of Algorithm \ref{algo:heuristic} for the data depicted in Figure \ref{fig:twogeneral}. The algorithm already terminates after the first iteration.} 
\label{fig:stages2}
\end{figure}

Figure \ref{fig:stages2} shows the first $\bar P_{\text{org}}$ and the result of the first run of Algorithm \ref{algo:nonmasssplit}. The approximate barycenter in the original support is a $1.102$--approximation of the exact barycenter, i.e., there is a $10.2\%$-error. This is improved to a $1.9\%$-error in the split-up using Algorithm \ref{algo:nonmasssplit}. The second iteration of Algorithm \ref{algo:heuristic} does not improve the solution anymore and the algorithm terminates.

This run completes in less than a second. In contrast, the computation of an exact barycenter for such a small problem size already is surprisingly hard: despite all measures having the same support, the set $S$ is of exponential size. The LP for an exact computation has $939510$ variables and $103032$ constraints and takes roughly $500$ seconds to solve \citep{abm-16}.

{\noindent \bf Observations.} Both of these sample runs are representative in a couple of ways. The approximation error for $\bar P_{\text{org}}$ and the first improvement to $\bar P'$ using Algorithm \ref{algo:nonmasssplit} are already much better than the guaranteed bound of $2$. In the computations in Section \ref{sec:scaling}, we have not encountered a run with an approximation error worse than $20\%$ for the initial $\bar P_{\text{org}}$ or $8.7\%$ for the initial improvement to $\bar P'$. The improvement between $\bar P_{\text{org}}$ and $\bar P'$ in the first iteration is significant. However, the additional iterations of Algorithm \ref{algo:heuristic} only perform minor improvements on the approximation factor. In fact, the example in Figure \ref{fig:stages} shows one of the largest improvements after the first iteration observed in all our computations.

Only two parts contribute significantly to the total running time of Algorithm \ref{algo:heuristic}: the first run of Algorithm \ref{algo:2approx} and the setup of the LP for Step $1$ in the second iteration. Together, these accounted for more than $80\%$ of the total running time. Only the first LP is run in full; the later LPs can be warm-started. Further, the setup of LPs in iterations $3$ or later is a simple update from the previous iteration. The split-up of mass in Algorithm \ref{algo:nonmasssplit} is negligible in running time.

An initial run of Algorithm \ref{algo:2approx} is unavoidable in all situations, but is efficient through the use of support set $S_{\text{org}}$. The setup of the second LP can be computationally expensive, because $S_{\text{org}}$ is replaced by the larger $\supp(\bar P')$. However, we observed that in practice $|\bar P'|$ does not only satisfy the guaranteed bound $|\bar P'|\leq (\sum_{i=1}^N |P_i|  - N + 1)^2$ (Theorem \ref{thm:fixingtheproperties}), but remains close to $\sum_{i=1}^N |P_i|  - N + 1$, the bound for $|P_{\text{org}}|$ (Theorem \ref{thm:heuristic}).

The main benefit of a full run of Algorithm \ref{algo:heuristic} is the recovery of the combination of a sparse support and a non-mass splitting transport. If these are not crucial properties for an application, we recommend performing just a single iteration of Algorithm \ref{algo:heuristic}, i.e., a single run of Algorithms \ref{algo:2approx} and \ref{algo:nonmasssplit}, for faster computations.


\subsection{Scaling of Algorithm \ref{algo:2approx} and runs for Algorithm \ref{algo:heuristic}}\label{sec:scaling}

Finally, we study the scaling of Algorithm \ref{algo:2approx} and the practical running time of our algorithms. Algorithm \ref{algo:2approx} is the main pillar of the two viable approaches in this paper for practical computations: either a single run of Algorithm \ref{algo:2approx} and \ref{algo:nonmasssplit} or a full run of Algorithm \ref{algo:heuristic}, if the problem size allows. It is the only LP without a warm-start. Algorithm \ref{algo:2approx} is based on an LP formulation using $S_{\text{org}}$ as the set of possible support points. Before we turn to more computations, let us take a closer look at the number of variables and constraints in this LP to set up proper expectations. Using $|P_{\text{sum}}|=\sum_{i=1}^N |P_i|$, this LP has 
$$|S_{\text{org}}|+|S_{\text{org}}|\cdot |P_{\text{sum}}| \text{ variables } \quad\quad \text{ and } \quad \quad N\cdot|S_{\text{org}}|+|P_{\text{sum}}| \text{ constraints.}$$
Two types of scaling are of interest: scaling the number $N$ of measures and scaling $|S_{\text{org}}|$. Note $N\cdot|S_{\text{org}}|\geq |P_{\text{sum}}|$, so the number of constraints scales linearly in $N$ and $|S_{\text{org}}|$. Further, note that the number of constraints is always lower than the number of variables ($|P_{\text{sum}}| \geq N$), and often dramatically so ($|P_{\text{sum}}| \gg N$). 

The dominating factor in the number of variables is $|S_{\text{org}}|\cdot |P_{\text{sum}}|$. If the support sets of the $P_i$ are disjoint, then $|S_{\text{org}}| =  |P_{\text{sum}}|$ and $|S_{\text{org}}|\cdot |P_{\text{sum}}|=|S_{\text{org}}|^2=|P_{\text{sum}}|^2$; the number of variables scales quadratically in $|S_{\text{org}}|$. For the sake of a simple analysis, we assume all measures have the same number of support points $|P_{\text{max}}|$. Then $|P_{\text{sum}}|=N\cdot |P_{\text{max}}|$ and the scaling of the number of variables is quadratic in $N$.


\subsubsection{Scaling for grid-structured data} 
First, we consider scaling for (two-dimensional) grid-structured data. Let $K$ denote the number of grid points in each direction. The grid has $K^2$ points, which is an upper bound on $|P_{\text{max}}|$ and on $|S_{\text{org}}|$. Thus the number of variables is bounded by $K^2 + K^2\cdot (N\cdot K^2)$ and the number of constraints is bounded by $N\cdot K^2 + N\cdot K^2= 2(N\cdot K^2)$. Note that the bound on the number of variables is roughly $K^4\cdot N$.

While the actual sizes of $|P_{\text{max}}|$ and $|S_{\text{org}}|$ are usually significantly smaller than $K^2$, they typically remain a linear fraction of $K^2$ (for MNIST digits between $\frac{1}{5}$ and $\frac{1}{3}$), and so these bounds immediately imply two types of consequences. First, the scaling with the number $N$ of measures is linear. Second, doubling the density $K$ of the underlying grids will multiply the number of variables by $16$.  

Let us turn to some computations. In Table \ref{table:avgruns}, we report on average errors and completion times for a large set of runs of Algorithm \ref{algo:heuristic}. For each of these runs, we used random samples of $16\times 16$ digits as the measures. Each row is based on the data from a total of $100$ runs, $10$ for each digit $0,1,\dots,9$. The table lists the initial error and time for the computation of an approximate barycenter $\bar P_{\text{org}}$ in $S_{\text{org}}$ (Algorithm \ref{algo:2approx}), the error for the first $\bar P'$ (Algorithm \ref{algo:nonmasssplit}), the error, time, and number of iterations for a full run of Algorithm \ref{algo:heuristic}, and the time for an exact barycenter computation.

We have been able to run Algorithm \ref{algo:heuristic} for up to $40$ measures in less than ten minutes, and Algorithm \ref{algo:2approx} for up to $100$ measures in less than twenty minutes. The big difference is the setup of the LP for the second iteration. In contrast, we have not been able to perform the computation of an exact barycenter for more than $8$ measures (within a fixed time limit of four hours), even using some refinements to an exact barycenter computation \citep{bs-18}. This is the reason for the "n{/}a" entries in the table, where an approximation error is not available because of the lack of an exact solution. For increasing number of measures, the difference between the running times of Algorithm \ref{algo:heuristic} and an exact computation becomes dramatic, even though grid-structured data, in fact, is a scenario where the exact LP does not scale exponentially (recall an exact barycenter is contained in an $N$-times finer grid).

The first row shows numbers on random samples of four measures, as in the example depicted in Figures \ref{fig:foursixes} and \ref{fig:stages}. We observed a termination of Algorithm \ref{algo:heuristic} after an average of $2.2$ iterations. This low number of iterations is not surprising because of the low initial error and our practical implementation of Step $3$ of Algorithm \ref{algo:nonmasssplit}; see Section \ref{sec:itersec}. The same effects extend to larger computations, where the approximation error of the initial $\bar P_{\text{org}}$ is already low, most of the further improvement already happens towards the first $\bar P'$, and less than $5$ iterations were necessary on average. The times include setup of the problems. We did not observe a clear pattern with respect to the errors for the first $\bar P_{\text{org}}$ and $\bar P'$ or the final approximation, but the average number of iterations of Algorithm \ref{algo:heuristic} increases slightly with the number of measures. 

\begin{table}[t]
\begin{center}
\begin{tabular}{c|c|c|c|c|c|c|c}
& \multicolumn{2}{|c|}{ first $\bar P_{\text{org}}$} & first $\bar P'$ &\multicolumn{3}{c|}{full run of Alg. \ref{algo:heuristic}} & \text{exact}\\ \hline
no. of measures & error & time (s) & error & error & time (s) & iterations & time (s) \\ \hline
4 & 14.8\% & 4.2 & 3.8\% & 3.1\% & 9.9 & 2.2 & 120 \\ \hline
5  &  15.2\%& 5.7 & 4.4\% & 4.1\%   & 16.4 & 2.8 & 204 \\ \hline
6 & 15.5\%& 8.5 &  4.5\% &  3.9\%  & 22.3 & 2.7 & 540\\ \hline
7 & 15.1\% & 12.1 & 4.6\% & 4.2\% & 29.8 & 3.1 & 1602 \\ \hline
8 & 16.2\% & 16.3 & 5.2\% & 4.8\%  & 36.7 & 3.1 & 4330\\ \hline
9 & n{/}a & 23.0 &  n{/}a & n{/}a & 45.2 & 3.4 & --\\ \hline
12 & n{/}a & 39.1 & n{/}a & n{/}a & 74.8 & 3.3 & --\\ \hline
16 & n{/}a & 58.4 & n{/}a & n{/}a & 99.3 & 3.7 & --\\ \hline
20 & n{/}a & 90.3 & n{/}a & n{/}a & 169.2 & 4.5 & --\\ \hline
40 & n{/}a & 298.7 & n{/}a & n{/}a & 557.3 & 5.0 & --\\ \hline
70 & n{/}a & 681.5 & n{/}a & -- & -- & -- & --\\ \hline
100 & n{/}a & 1198.2 & n{/}a & -- & -- & -- & --\\ \hline
\end{tabular}
\end{center}\caption{Average numbers (error, time) for an initial approximation $\bar P_{\text{org}}$, first $\bar P'$, full runs of Algorithm \ref{algo:heuristic}, and an exact computation in a grid setting. The numbers in each row were derived from $100$ random samples of $16\times 16$ digits from the MNIST data set.} \label{table:avgruns}
\end{table}

\begin{figure}[h]
\hspace*{-0.5cm}
\subfloat[$\bar P_{\text{org}}$, $32\times32$ digits]{
  \includegraphics[scale=0.28]{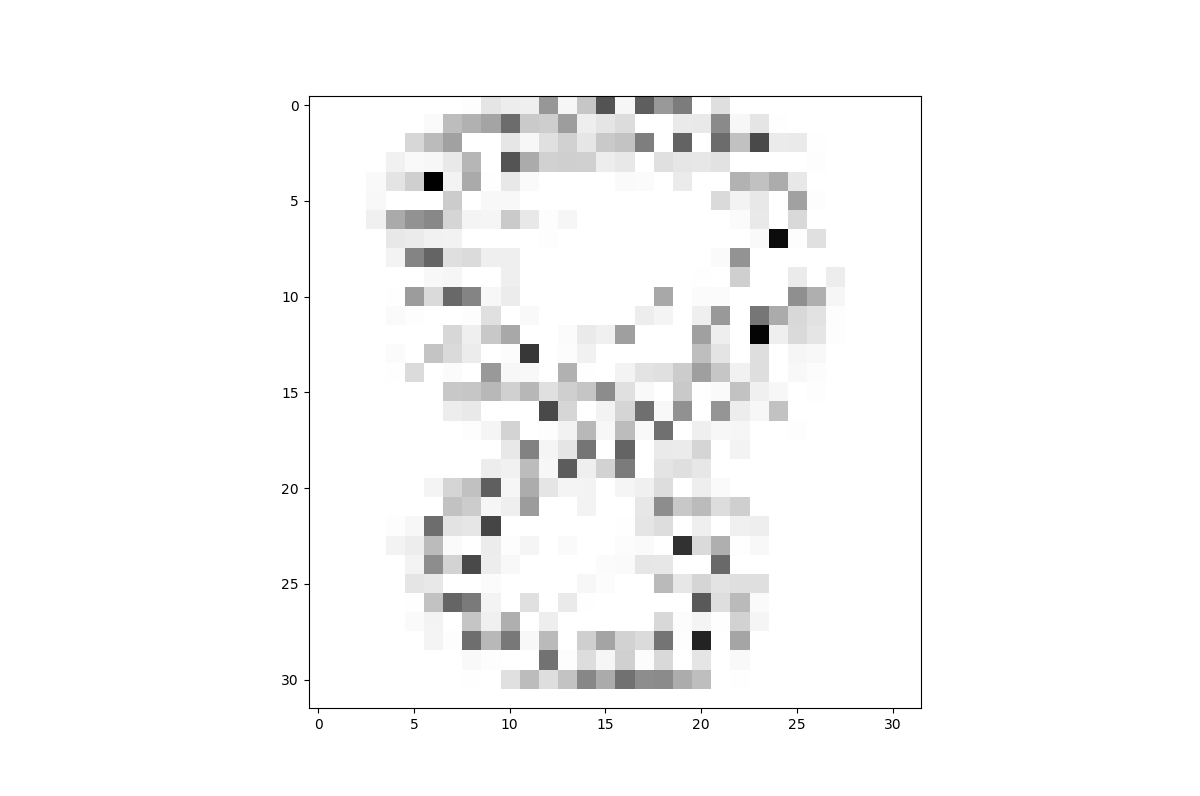}
}\hspace*{-2cm}
\subfloat[$\bar P_{\text{org}}$, $64\times64$ digits]{
  \includegraphics[scale=0.28]{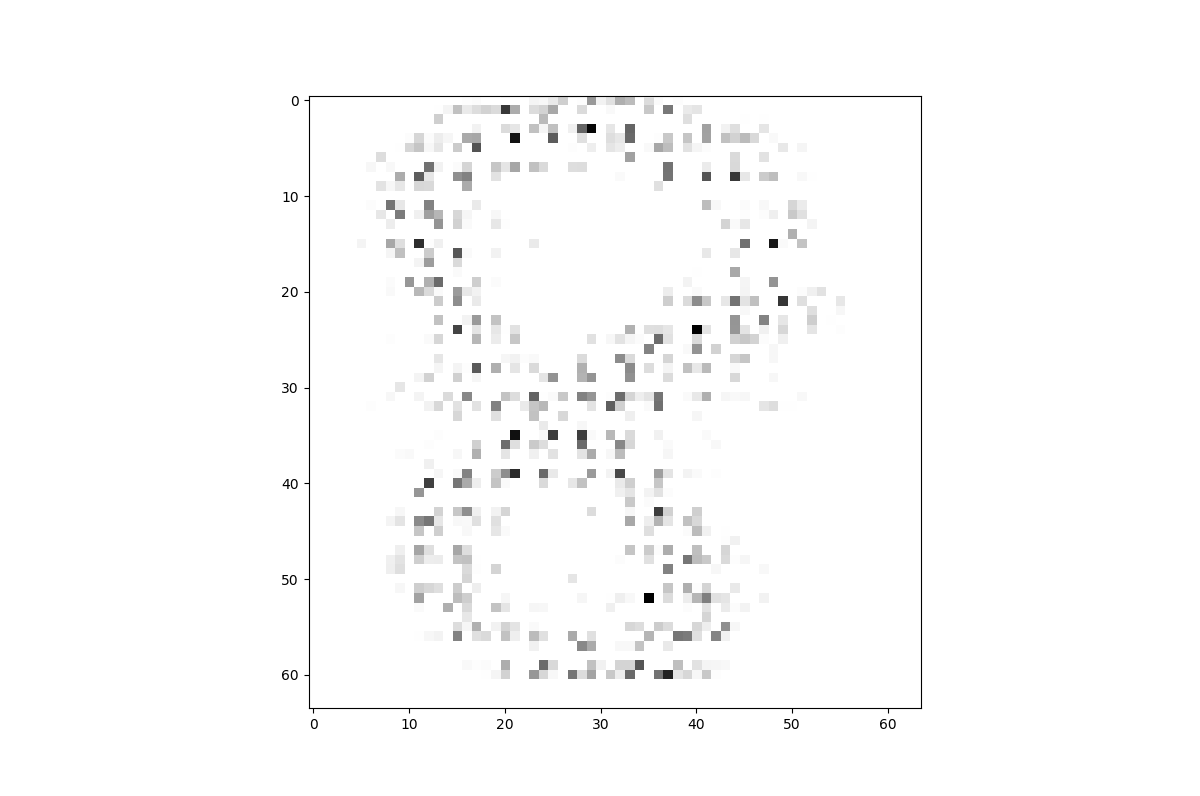}
}

\caption{Approximate barycenters $\bar P_{\text{org}}$ for a run of Algorithm \ref{algo:2approx} for $4$ digits in $32\times 32$ and $64\times 64$ grids. These computations already took several minutes, respectively more than an hour. Computations in denser grids quickly become impossible due to the quadratic scaling of the underlying LP with respect to $|S_{\text{org}}|$.} 
\label{fig:largergrids}
\end{figure}

Computations in denser grids quickly become impossible for the algorithms in this paper, due to the quadratic scaling of the underlying LP with respect to $|S_{\text{org}}|$. Recall that doubling the density $K$ of a grid multiplies the number of variables by $16$. Figure \ref{fig:largergrids} shows the results of Algorithm \ref{algo:2approx} for four measures in a $32\times32$ grid and a $64\times 64$ grid. The computations took about $5$ minutes, respectively $92$ minutes. The $64\times 64$ example exceeds $10$ million variables and is only solvable because of the extremely low number of constraints.

\begin{table}[t]
\begin{center}
\begin{tabular}{c|c|c|c|c|c|c|c}
& \multicolumn{2}{|c|}{ first $\bar P_{\text{org}}$} & first $\bar P'$ &\multicolumn{3}{c|}{regularization-based  $\bar Q$}\\ \hline
no. of measures & error & time (s) & error & error & time (s) & error gap {\small$\phi(\bar Q)\slash\phi(\bar P)$} \\ \hline
4 & 14.8\% & 4.2 & 3.8\% & 22.4\% & 6.8 & 1.179 \\ \hline
6 & 15.5\%& 8.5 &  4.5\% &  21.9\%  & 10.2 & 1.167 \\ \hline
8 & 16.2\% & 16.3 & 5.2\% & 23.6\%  & 13.7 & 1.175 \\ \hline
12 & n{/}a & 39.1 & n{/}a & n{/}a & 20.4 & 1.182 \\ \hline
16 & n{/}a & 58.4 & n{/}a & n{/}a & 27.7 & 1.188 \\ \hline
20 & n{/}a & 90.3 & n{/}a & n{/}a & 33.0 & 1.179 \\ \hline
40 & n{/}a & 298.7 & n{/}a & n{/}a & 67.9 & 1.191 \\ \hline
70 & n{/}a & 681.5 & n{/}a &  n{/}a & 115.5 & 1.185 \\ \hline
100 & n{/}a & 1198.2 & n{/}a &  n{/}a & 164.9 & 1.182 \\ \hline
\end{tabular}
\end{center}\caption{A comparison (error, time) of an initial approximation $\bar P_{\text{org}}$ and first $\bar P'$ to the result $\bar Q$ of a widely-used regularization-based approximation algorithm (\cite{cd-14}). The error gap shows the ratio $\phi(\bar Q)\slash\phi(\bar P)$ of objective function values for $\bar Q$ and $\bar{P'}$.} \label{table:comparison} 
\end{table}

For grid-structured data, there are many algorithms in the literature that are much faster options to tackle larger problem instances \citep{brpp-15,c-13,cd-14,pc-19}. In Table \ref{table:comparison}, we compare approximation errors and computational speed for the first $\bar P_{\text{org}}$ and $\bar P'$ to the widely-used, regularization-based algorithm from \cite{cd-14}. The algorithm uses Sinkhorn distances to simplify the objective function and leads to dense approximations $\bar Q$ over $S_\text{org}$. Recall that $\bar P_{\text{org}}$ is already an exact, sparse solution over $S_\text{org}$, and is further refined to $\bar P'$. Thus, the approximation error for $\bar Q$ always has to be worse than for $\bar P_{\text{org}}$, and the gap becomes signficant when compared to $\bar P'$. 

For regularization-based algorithms, there are several parameters that allow for a tradeoff of computational speed and quality of result. We used the recommended settings for MNIST data from \cite{cd-14}. Table \ref{table:comparison} shows exact approximation errors for the small instances for which an exact barycenter computation  was possible (up to $8$ measures), as well as the ratio $\phi(\bar Q)\slash\phi(\bar P)$ between cost of transport for $\bar Q$ and $\bar P$ for all instances. Our computations revealed a ratio of about $1.18$, respectively a gap of $18\%$, throughout. In view of computational speed, however, the regularization-based algorithm scales dramatically better with the number of measures than our approach, as expected. We observe a near-linear increase of running time, which sharply contrasts with a linear scaling of the size of the LP for Algorithm \ref{algo:2approx}.

Of course, one has to be careful in this comparison. In addition to the better approximation error, the additional computational cost of the algorithms in this paper leads to several favorable properties that are hard to measure quantitatively: the guarantee of a $2$-approximation, sparsity, non-mass split, numerical stability, and support in $S$ and not only in $S_\text{org}$. For grid-structured data, the tradeoff to obtain these properties may often not be worth the additional time in practice. However, our algorithms run without specification of a fixed support set for the solution, and thus have the ability to work for any data.  Next, we turn to a best-case type of data for our algorithms, where we are able to scale our computations to thousands of measures. Notably, it is data in general position, for which many algorithms in the literature do not work at all.

\subsubsection{Scaling for data in general position} 
Applications in operations research often are based on a small set of geographical locations that do not exhibit an obvious structure. In this section, we consider data sets as depicted in Figure \ref{fig:twogeneral} - a set of $N$ measures that all have the {\em same}, small support of size $|P_{\text{max}}|$. This configuration leads to a best-case scenario in that $|S_{\text{org}}|=|P_{\text{max}}|$, i.e., the union of original supports is the same as any of the supports. This has an extremely positive effect on the size of the LP for Algorithm \ref{algo:2approx}, respectively the first iteration of Algorithm \ref{algo:heuristic}: the number of variables is $|P_{\text{max}}|+N\cdot |P_{\text{max}}|^2$ and the number of constraints is $2(N\cdot |P_{\text{max}}|)$. In comparison, if the support points of the $P_i$ did not overlap, one would have $N\cdot |P_{\text{max}}|+(N\cdot |P_{\text{max}}|)^2$ variables and  $N\cdot |P_{\text{max}}|+N^2\cdot |P_{\text{max}}|^2$ constraints. The advantage is a factor of about $N$ in the number of both variables and constraints. Note, however, that the same does not transfer to the LPs in later iterations of Algorithm \ref{algo:heuristic}, where $\text{supp}(\bar P')$ is used as the new support, and is only guaranteed to satisfy the bound in Theorem \ref{thm:fixingtheproperties}.

\begin{figure}
\begin{center}
\subfloat[$9$ support points, $48$ measures]{
  \includegraphics[scale=0.76]{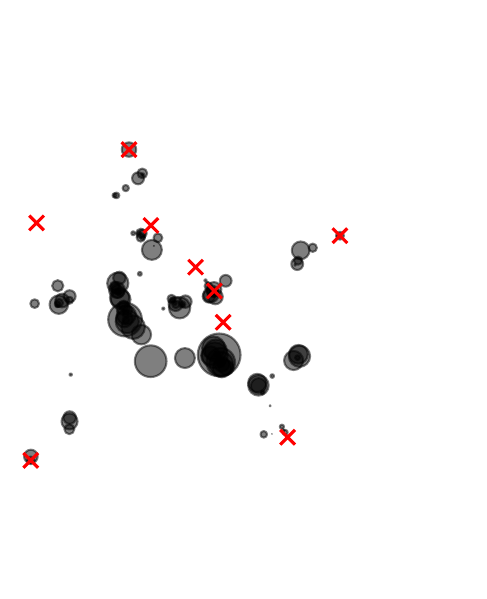}
}
\subfloat[$9$ support points, $96$ measures]{
  \includegraphics[scale=0.76]{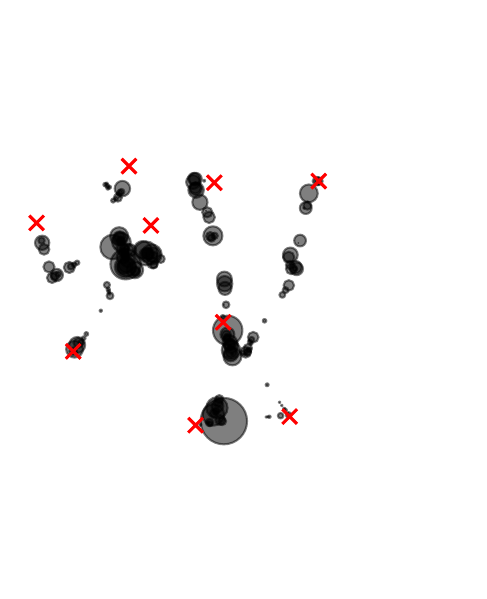}
}

\subfloat[$13$ support points, $36$ measures]{
  \includegraphics[scale=0.76]{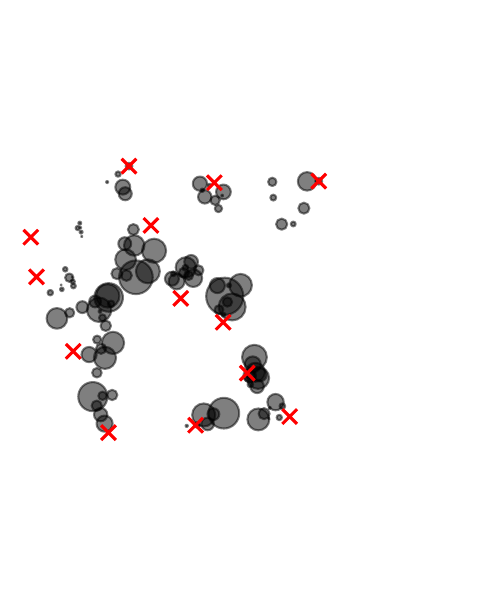}
}
\subfloat[$17$ support points, $36$ measures]{
  \includegraphics[scale=0.76]{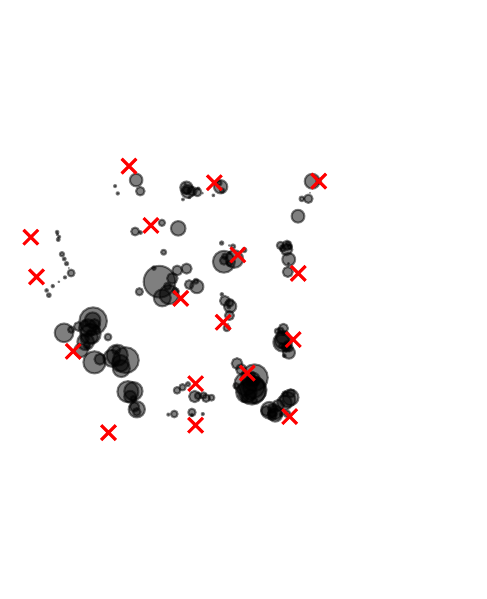}
}
  \end{center}
\caption{Approximate barycenters $\bar P'$ at the end of the first iteration of Algorithm \ref{algo:heuristic} for different support sets and number of measures. Full runs of Algorithm \ref{algo:heuristic} completed in less than $30$ seconds.} 
\label{fig:biggerexamples}
\end{figure}

In Table \ref{table:avgruns2}, we report on average errors and completion times for runs of Algorithm \ref{algo:2approx} and Algorithm \ref{algo:heuristic}. The numbers in each row were derived from $100$ runs for the given number of measures, all consisting of the same $|P_{\text{max}}|=9$ support points in general position, randomly generated for each run. (We discuss the effect of scaling to larger $|P_{\text{max}}|$ below; here we exhibit a best-case scenario.) The measures were constructed through a random assignment of varying masses to the support. We also chose the weights $\lambda_i$ randomly.  Figure \ref{fig:biggerexamples} (top) shows some approximate barycenters computed in this setting.

The combination of data in general position, which makes a discretization of the underlying space unavailable, and a non-uniform weight vector makes for an impractical setting for algorithms in the literature; see Section \ref{sec:exact}. This is why our comparisons are restricted to an exact, LP-based solution. However, exact computations in this setting are extremely hard due to the (always) exponential scaling of $S$ and the corresponding LPs. (Here, $|S|=|P_{\text{max}}|^N$.) The largest number of measures for which we successfully found an exact solution is $12$ (in a bit less than 3.5 hours). Because of this, approximation errors are not available for more than $12$ measures. At the same time, the LP for Algorithm \ref{algo:2approx} for such an instance still is of trivial size: it has $12\cdot 9^2$ variables and $2(12\cdot 9)$ constraints. The speed-up over an exact computation for a set of $8$ measures (the example in Section \ref{sec:scaling}) is a factor of more than $600$. This factor escalates quickly - for $12$ measures our algorithm is already faster by a factor of more than $4000$. 

\begin{table}[h]
\begin{center}
\begin{tabular}{c|c|c|c|c|c|c|c}
 & \multicolumn{2}{|c|}{first $\bar P_{\text{org}}$} &  first $\bar P'$& \multicolumn{3}{c|}{full run of Alg. \ref{algo:heuristic}} & \text{exact}\\ \hline
no. of measures & error & time (s) & error & error & time (s) & iterations & time (s) \\ \hline
8 & 10.1\% & 0.7 & 2.0\% & 1.6\% & 0.9 & 1.4 & 505\\ \hline
12  &  9.8\%& 1.1 & 2.4\% & 1.9\%  & 3.0 & 1.8 & 12400\\ \hline
50 & n{/}a & 2.3 &  n{/}a & n{/}a & 8.0 & 2.2 & --\\ \hline
100  & n{/}a & 3.5 &  n{/}a & n{/}a & 26.6 & 2.3 & --\\ \hline
200  & n{/}a & 5.3 &  n{/}a & n{/}a & 112.8 & 2.9 & --\\ \hline
500  & n{/}a & 7.9 &  n{/}a & n{/}a & 730.4 & 3.5 & --\\ \hline
1000  & n{/}a & 18.5 &  n{/}a & --  & -- & -- & -- \\ \hline
5000  & n{/}a & 92.0 &  n{/}a & --  & -- & -- & -- \\ \hline
10000  & n{/}a & 229.2 &  n{/}a & --  & -- & -- & -- \\ \hline
20000  & n{/}a & 963.4 &  n{/}a & --  & -- & -- & -- \\ \hline
\end{tabular}
\end{center}\caption{Average numbers (error, time) for an initial approximation $\bar P_{\text{org}}$, first $\bar P'$, and full runs of Algorithm \ref{algo:heuristic}, and an exact computation for data in general position. The numbers in each row were derived from $100$ runs on a set of measures with $9$ support points of randomly chosen masses.} \label{table:avgruns2}
\end{table}

We have been able to run Algorithm \ref{algo:2approx} for up to $20000$ measures in about 15 minutes. Instances up to $1000$ measures solve in less than $20$ seconds. One of the main reasons for these low running times is the extremely low number of constraints. The same scalability of Algorithm \ref{algo:heuristic} cannot be expected, as the size of the support for LPs in later iterations is equal to $|\bar P'|$, which is between linear and quadratic in $N\cdot |P_{\text{max}}|$  (Corollary \ref{sparseequ2}, Theorem \ref{thm:fixingtheproperties}). We were able to run it for up to $500$ measures in a bit more than $12$ minutes. Scaling further, one quickly reaches a point where the LP for the second iteration cannot be constructed anymore. (For $1000$ measures, it would have close to $100$ million variables.) 

 The times reported in Table \ref{table:avgruns2} include the setup of the problems. Again, we observe that $\bar P_{\text{org}}$ is significantly improved in the first step towards $\bar P'$, and that further iterations do not change it noticably anymore. The average number of iterations increases slowly with the number of measures. The increase is slower than in the grid-structured setting. Informally, less repetition in the weighted centroids for different combinations of support points means that later iterations of Algorithm \ref{algo:heuristic} are less likely to further improve the solution.

Unlike for grid-structured data, a scaling of $|P_{\text{max}}|$ for Algorithm \ref{algo:2approx} is easier in this setting. Figure \ref{fig:biggerexamples} (bottom) show two examples with a larger support. Recall that the number of variables for Algorithm \ref{algo:2approx} is quadratic in $|P_{\text{max}}|$. Doubling the number of support points increases the number of variables in the problem by factor $4$, the same effect as increasing the number of measures by factor $4$. The number of constraints is doubled (and remains extremely low in comparison to the number of variables). An instance with $2500$ measures of $18$ support points exhibits a similar running time to an instance with $10000$ measures of $9$ support points (line 9 in Table \ref{table:avgruns2}).

Summing up, the combination of few, overlapping support points in general position and a large number of measures is a best-case scenario for the combination of a single run of Algorithms \ref{algo:2approx} and \ref{algo:nonmasssplit}. Applications in operations research, such as facility location problems, often fall into this category. In this setting, exact computations are impossible for any reasonable problem size. Algorithms in the literature are not designed to work (or work well) for such data. In this situation, we recommend use of the presented methods for computational speed. In all other situations, the favorable properties of the output - like sparsity, non-mass split, and a guaranteed error bound - have to be crucial to the application to be worth the significant additional computational cost over popular heuristics.

{\small \section*{Acknowledgments} 

\noindent The author would like to thank Ethan Anderes for the support with implementations in the Julia language, and Jacob Miller for the helpful discussions. The author gratefully acknowledges support through the Collaboration Grant for Mathematicians  {\em Polyhedral Theory in Data Analytics} of the Simons Foundation.}


\bibliography{barycenters_literature}
\bibliographystyle{plain}

\section*{Appendix}

\appendix
\normalsize
\section{- Proofs of Theorems \ref{thm:fixingtheproperties} and \ref{cor:nonmasssplit}}\label{sec:localimprovement}

We begin by proving Theorem \ref{thm:fixingtheproperties}.

\setcounter{theorem*}{2}
\begin{theorem*}
Algorithm \ref{algo:nonmasssplit} returns a measure $\bar P'$ supported on a subset of $S$ with $\phi(\bar P') \leq 2 \cdot \phi(\bar P)$ and there is a non-mass splitting transport realizing this bound. Further  $|\bar P'|\leq (\sum_{i=1}^N |P_i|  - N + 1)^2$.
\end{theorem*}

\proof{}
First, note that the $P_i^l$ constructed in Step $1$ satisfy $\supp(P_i^l)\subset \supp(P_i)$.
Thus $\supp(\bar P^l)\subset S$, and consequently $\supp(\bar P')\subset S$. Further, $\bar P'= \sum_{l=1}^r \bar{P^l}$ is a measure. This holds because $\sum_{l=1}^{r} d_l= \sum_{l=1}^{r} z_{t_l} = 1$, because Step $2$ does not affect this sum, and because the total mass in $\bar{P^l}$ equals $d_l$ by construction. Thus, $\bar P'$ is a measure supported in $S$.

Second, we prove correctness of Step $2$. We will show that a greedily lexicographically maximal $(d_1,\dots,d_r)$ is created while retaining an approximate barycenter in $\supp(P_{\text{org}})$. In particular, we have to show that the objective function value $\phi(\bar P_{\text{org}})$ does not change during the shift of mass. For a simple wording, let $\bar P_\text{lex}$ be the measure corresponding to $(d_1,\dots,d_r)$ after Step $2$. We will prove $\phi(\bar P_{\text{org}})=\phi(\bar P_{\text{lex}})$.

Let $x_{iq_i}^l\in P_i^l$ for $i\leq N$ and $c=\sum_{i=1}^N \lambda_i x_{iq_i}^l$, as in Step $2 a)$. Then $\|c-s_l\|\leq \|c-s_j\|$ for all $j\neq l$. To see this, recall $$\sum\limits_{i=1}^N \lambda_i \|s-x^l_{iq_i}\|^2 = \sum\limits_{i=1}^N\lambda_i (\|s-c\|^2+\|c-x^l_{iq_i}\|^2),$$
as demonstrated in the proof of Theorem \ref{thm:originalsupport}. If $\|c-s_l\|> \|c-s_j\|$ for some $j\neq l$, $\bar P_\text{org}$ would not have been optimal. 

By $q_i= \text{arg}\max_{q\leq |P_i^l|} (s_j-s_l)^Tx^l_{iq}$ in Step $2a)$, we pick the $x^l_{iq_i}$ such that their weighted centroid $c=\sum_{i=1}^N \lambda_i x^l_{iq_i}$ maximizes the difference $\|c-s_l\|^2 - \|c-s_j\|^2 \leq 0$. Only if  $\|c-s_l\|^2= \|c-s_j\|^2$, mass is shifted from $s_l$ to $s_j$. But then the approximation error does not change, because
 $$\sum_{i=1}^N \lambda_i \|s_j-x^l_{iq_i}\|^2  = \sum\limits_{i=1}^N\lambda_i (\|s_j-c\|^2+\|c-x^l_{iq_i}\|^2) = \sum\limits_{i=1}^N \lambda_i \|s_l-x^l_{iq_i}\|^2 .$$
Thus, the objective function value does not change during Step $2$; we have $\phi(\bar P_{\text{org}})=\phi(\bar P_{\text{lex}})$.

By definition of the running indices $l$ and $j$, mass can only be moved from support points of higher index $l$ to support points of lower index $i$. For each pair of $l$ and $j$, we repeat this shift of mass until there is no weighted centroid with $\|c-s_l\|= \|c-s_j\|$ anymore. Due to decreasing $l$ in the outer loop and increasing $j$ in the inner loop, $(d_1,\dots,d_r)$ is transformed to be greedily lexicographically maximal and the corresponding measure remains an approximate barycenter.

Next, we prove correctness of Steps $3$ and $4$. We show that $\phi(\bar P_{\text{org}})\geq\phi(\bar P')$. Further,  we show that for each constructed partial measure $\bar P^l$ there is a non-mass splitting transport to the $P_i^l$, and that they combine to a $\bar P'$ that allows for a non-mass splitting transport that is at least as good as an optimal transport for $\bar P_{\text{org}}$. Finally, we show $|\bar P'|\leq (\sum_{i=1}^N |P_i|  - N + 1)^2$. 

Recall that in Step $3$, the mass of each $s_l$ is spread out to a set of weighted centroids to obtain $\bar P^l$. Independently of how the $x^l_{iq_i}$ are picked from the $P^l_i$ for all for all $i\leq N$, their weighted centroid $c=\sum_{i=1}^N \lambda_i x^l_{iq_i}$ satisfies $\sum_{i=1}^N \lambda_i \|c-x_{iq_i}^l\|^2 \leq \sum_{i=1}^N \lambda_i  \|s_l-x_{iq_i}^l\|^2$. By construction of $\bar P'$ from the $\bar P^l$ (Step $4$), this already implies $\phi(\bar P') \leq \phi(\bar P_{\text{org}})$. The algorithm started with a $2$-approximation, and thus it is guaranteed to return a $\bar P'$ with $\phi(\bar P') \leq 2 \cdot \phi(\bar P)$.

The existence of a non-mass splitting transport from $\bar P'$ to $P_1,\dots,P_N$, and the fact that this transport realizes the above bound, is a consequence of two reasons. First, each $\bar P^l$ itself allows for a non-mass splitting transport to the $P_i^l$ by lexicographically maximal choice of the $x_{iq_i}^l$ in Step $3a)$: due to this choice, the first constructed weighted centroid $c$ is lexicographically maximal among all (possible) weighted centroids that can be constructed from any $x_{iq}^l$ in the $P^l_i$. Further, by reducing the mass at each used support point by $d_{\text{min}}$ in Step $3b)$, at least one of the $d^l_{iq_i}$ becomes $0$. The corresponding support point is removed from $P^l_i$ (followed by some reindexing) and thus cannot be used for the construction of a weighted centroid in further iterations. Thus, the second centroid constructed in the inner loop is lexicographically {\em strictly} smaller than the first one. The same holds for all subsequent ones. 

Second, any two partial measures $\bar P^{l_1}$,  $\bar P^{l_2}$ from Step $3$ satisfy $\text{supp}(\bar P^{l_1}) \cap \text{supp}(\bar P^{l_2}) = \emptyset$ for $l_1\neq l_2$, because of the earlier preprocessing in Step $2$: weighted centroids that would be equally distant from both $s_{l_1}$ and $s_{l_2}$ cannot exist, because this would have caused a shift of mass to the lower index in Step $2$ to create a lexicographically larger $(d_1,\dots,d_r)$. Summing up, $\bar P'$ consists of a set of distinct support points, for which it is trivial to give a non-mass splitting transport to the $P_i$ that is at least as good as an optimal transport for $\bar P_{\text{org}}$: this transport just sends the whole mass of each support point in $\bar P'$ to the support points in the $P_i$ that were used for its construction.

The removal of at least one support point from a $P^l_i$ in Step $3b)$ implies that there are at most $\sum_{i=1}^N |P^l_i| - N + 1$ runs of $3a)$ and $3b)$ to construct a  $P^l$: the 'go back to $a)$' statement is applied while $d_l>0$; this is the case while there still is a support point in a $P^l_i$ with mass on it. In the final run of Steps $3a)$ and $3b)$ for each $P^l$, all the $P^l_i$ have precisely one support point with the same mass left. This gives the claimed bound, and in particular $|P^l|\leq \sum_{i=1}^N |P^l_i| - N + 1$. 

Due to $|P^l_i|\leq |P_i|$ and $|\bar P_\text{org}| \leq \sum_{i=1}^N |P_i| -N  + 1$, we obtain $$|\bar P'|=\sum\limits_{l=1}^{|\bar P_\text{org}|}  |\bar P^l| \leq \sum\limits_{l=1}^{|\bar P_\text{org}|} (\sum\limits_{i=1}^N |P^l_i| -N  + 1) \leq \sum\limits_{l=1}^{|\bar P_\text{org}|} (\sum\limits_{i=1}^N |P_i| -N  + 1) \leq (\sum_{i=1}^N |P_i|  - N + 1)^2.$$ Thus $\bar P'$ satisfies all claimed properties. \hfill$\square$ 
\endproof

\noindent Next, we prove that Algorithm \ref{algo:nonmasssplit} runs in strongly-polynomial time.

\begin{theorem*}
For all rational input, a measure can be computed in strongly-polynomial time that is a $2$-approximation of a barycenter and for which there is a non-mass splitting transport realizing this bound. 
\end{theorem*}

\proof{}
We consider the running time of each part of the algorithm. For readability, we say `polynomial' in this proof in place of `strongly-polynomial'. We use `linear' and `quadratic' to refer to the bit size $\mathcal{I}$ of the input. Note that $N$, the $|P_i|$, and the dimension $d$ are all bounded above by $|\mathcal{I}|$.

 In Step $1$, the input for the subsequent steps is created. By sparsity of $\bar P_\text{org}$, $r \leq \sum_{i=1}^N |P_i| - N + 1$. For each of the $r$ support points $s_l$, $N$ images $P^l_i$ with $|P^l_i|\leq |P_i|$ are created. In the application of the stated rule, each $y_{it_lk}$ has to be processed (at most) once. For each $y_{it_lk}$, a single comparison and a fixed number of elementary operations suffices to update the support point and mass in $P_i^l$. In total, data structures of polynomial size are created in polynomial time.

Step $2$ is the preprocessing of $(d_1,\dots,d_r)$ to be greedily lexicographically maximal. For each pair of support points $s_l, s_j$ with $j < l$, we perform the inner part of the loop. Finding $q_i$ in $2a)$ can be done by considering all $x^l_{iq} \in P_i^l$ exactly once and comparing the inner products $(s_j-s_l)^Tx^l_{iq}$. This is possible in linear time. $c$ is created through the scaling and the sum of $N$ rational $d$-dimensional vectors. 

Step $2b)$ begins with the computation of $c-s_j$ and $c-s_l$, then computes $\|c-s_j\|^2=(c-s_j)^T(c-s_j)$ and  $\|c-s_l\|^2=(c-s_l)^T(c-s_l)$, and then compares the two values. This is possible in quadratic time. Picking the minimal mass among the $x_{iq_i}^l$ is possible in linear time, and so is updating the masses, performing the set operations on $P_i^l$ and $P_i^j$, and reindexing. In this update, $|P^l_i|$ is reduced by at least one, so the 'go back to $a)$' statement is followed at most $|P^l_i|$ times. Summing up, Step $2$ runs in polynomial time.

Step $3$ performs the spread-out of the $r$ support points. Picking a lexicographically maximal support point $x_{iq_i}^l$ in $3a)$ can be done by considering all support points in $P_i^l$ once. One saves the current best support point and compares each other support point with respect to their lexicographic order. For identifying the lexicographic order of a pair of $d$-dimensional support points, (at most) all $d$ of their coefficients have to be compared to each other. This is possible in linear time. Again, $c$ is created through the scaling and the sum of $N$ rational $d$-dimensional vectors. 

In $3b)$, we pick the minimal mass among the $x_{iq_i}^l$ used for the construction of $c$, which can be done in linear time. The same holds for the update of masses, the set operations on $P_i^l$, and the reindexing. By this update, the size of one of the $|P^l_i|$ is reduced by at least one, so the 'go back to $a)$' statement is followed not more than $\sum_{i=1}^N |P^l_i|$ times; more precisely, there are at most $|P^l_i| - N + 1$ runs of $3a)$ and $3b)$ for each $l$. Summing up, the construction of each $\bar P^l$ runs in polynomial time, and so does the construction of all the $\bar P^l$.  

In Step $4$, the partial measures $\bar P^l$ are combined to obtain $\bar P'$. This is the construction of a measure with the appropriate mass put on at most $|\bar P'|\leq (\sum_{i=1}^N |P_i|  - N + 1)^2$ support points. Each of these support points is just a copy of a support point in one of the $\bar P^l$. Thus, all steps run in polynomial time, which proves the claim. \hfill$\square$ 
\endproof

\section{- Proof of Theorem \ref{thm:heuristic}}\label{sec:iterativeproofs}

\begin{theorem*}
Algorithm \ref{algo:heuristic} returns an approximate barycenter $\bar P'$ supported on a subset of $S$ for which $\phi(\bar P') \leq 2 \cdot \phi(\bar P)$, where  $\bar P$ is a barycenter, and there is a non-mass splitting optimal transport realizing this bound. Further  $|\bar P'|\leq \sum_{i=1}^N |P_i|  - N + 1$.
\end{theorem*}

\proof{}
First, recall that the output $\bar P'$ of Algorithm \ref{algo:nonmasssplit} (Step $2$) always satisfies $\text{supp}(\bar P') \subset S$. Further, Algorithm \ref{algo:nonmasssplit} always returns a measure that has a corresponding non-mass splitting transport. As $\bar P_\text{org}$ from Algorithm \ref{algo:2approx} (Step $1$) is not changed in the final run of Algorithm \ref{algo:nonmasssplit}, the returned non-mass splitting transport is optimal. Further, recall that all approximate barycenters $\bar P_\text{org}$ computed in Step $1$ have a support that satisfies $|\bar P_\text{org}|\leq \sum_{i=1}^N |P_i|  - N + 1$. This transfers to the sparsity of $\bar P'$ returned by Algorithm \ref{algo:heuristic}.

It remains to prove termination of Algorithm \ref{algo:heuristic} and the error bound. We will do so by showing that $\phi(\bar P') < \phi(\bar P_\text{org})$ if $\bar P' \neq \bar P_\text{org}$ for $\bar P_\text{org},\bar P' $ from the same iteration. This leads to a strictly decreasing sequence of values $\phi(\bar P')$ as long as the algorithm keeps running. The first approximate barycenter in this sequence already is a $2$-approximation and it can only become better. This immediately gives $\phi(\bar P') \leq 2 \cdot \phi(\bar P)$. At the end of each Step $2$, we update $S_\text{org}=\text{supp}(\bar P')\subset S$ before going back to Step $1$, where an exact optimum over this new support, a subset of $S$, is computed. Because of this, and the fact that there are only finitely many subsets of $S$, the sequence of values $\phi(\bar P')$ is finite.

Now, it only remains to prove that $\phi(\bar P') < \phi(\bar P_\text{org})$ if $\bar P' \neq \bar P_\text{org}$. We begin by considering Step $3$ of Algorithm \ref{algo:nonmasssplit}. Assume $P^l_i$ consists of a single support point $x^l_{i1}$ for all $i\leq N$. Then the unique barycenter $\bar P^l$ of the $P^l_i$ is the weighted centroid $c = \sum_{i=1}^N \lambda_i x^l_{i1}$ and the cost of transport from $\bar P^l$ to all the $P^l_i$ is $\phi(\bar P^l) = d_l\cdot \sum_{i=1}^N\lambda_i \|c-x^l_{i1}\|^2$. For all $s\neq c$, in particular for $s=s_l$, we get
$$\phi(\bar P^l) = d_l\cdot \sum\limits_{i=1}^N\lambda_i \|c-x^l_{i1}\|^2 < d_l\cdot \sum\limits_{i=1}^N\lambda_i \|s-x^l_{i1}\|^2.$$

If some of the $P^l_i$ consist of more than one support point, Step $3$ selects a set of exactly one support point $x_{iq}^l$ from each measure $P^l_i$, forms a weighted centroid $c$ with corresponding mass $d_c=d_{\text{min}}$, and adds it to $\supp(\bar P^l)$. Then this scheme is repeated for the remaining support points and remaining mass. This means that $\bar P^l$ is constructed as a set of weighted centroids $c$ of support points $x^l_{iq}$ to which these centroids $c$ transport. Each of them satisfies  $d_c\cdot \sum_{i=1}^N\lambda_i \|c-x^l_{iq}\|^2 \leq d_c\cdot \sum_{i=1}^N\lambda_i \|s_l-x^l_{iq}\|^2$. By summing over all $c$ that are constructed, one obtains $$\phi(\bar P^l)\leq \sum_{i=1}^N\lambda_i \sum_{q=1}^{|P_i^l|} d^l_{iq} \cdot \|s_l-x_{iq}^l\|^2.$$

  Informally, it is at least as costly to transport to the measures $P^l_i$ from the support point  $s_l$ as from the set of weighted centroids (with appropriate masses) constituting $\bar P^l$. Equality in the above can only hold if the single support point $s_l$ itself already is the weighted centroid of single-support point measures $P^l_1,\dots,P^l_N$. But this means that Step $3$ of Algorithm $2$ just copies $s_l$ with mass $d_l$ to $\bar P^l$. The algorithm stops when  $\bar P'=\bar P_\text{org}$. By $\phi(\bar P')= \sum_{l=1}^r  \phi(\bar P^l)$, this means all $s_l$ have to satisfy $\phi(\bar P^l)= \sum_{i=1}^N\lambda_i \sum_{q=1}^{|P_i^l|} d^l_{iq} \cdot \|s_l-x_{iq}^l\|^2.$ So all $s_l$ are already the weighted centroids of their single-support measures $P^l_i$. 

Further, note that when a shift of mass from $s_l$ to $s_j$ with $j<l$ happens in Step $2$ of Algorithm \ref{algo:nonmasssplit}, then Step $3$ is guaranteed to find a strictly better transport than before: there exists a set of support points that, before the shift, receive transport from $s_l$, but have a weighted centroid $c\neq s_l$. Such a set of support points would be moved from $P^l_i$ to $P^j_i$ (and at least one of the support points was not associated to $s_j$ before). Then $s_j$ is guaranteed to split mass and, in the following Step $3$, the cost of transport is strictly improved; see above.

Thus $\phi(\bar P') < \phi(\bar P_\text{org})$ if $\bar P' \neq \bar P_\text{org}$ and Algorithm \ref{algo:heuristic} terminates with $\bar P' = \bar P_\text{org}$ in the final iteration. \hfill$\square$ 
\endproof

\end{document}